\documentclass[10pt, francais]{smfart}
\usepackage{amsfonts, amssymb, amsmath, amsthm, latexsym, enumerate, epsfig, color, mathrsfs, fancyhdr}
\usepackage[all]{xy}

\usepackage[frenchb]{babel}

\usepackage{mysections}

\usepackage[mac]{inputenc}

\usepackage{hyperref}
\hypersetup{bookmarksnumbered,bookmarksopen,bookmarksopenlevel=2}

\def\N{\ensuremath{\mathbf{N}}}

\def\Q{\ensuremath{\mathbf{Q}}}
\def\R{\ensuremath{\mathbf{R}}}
\def\C{\ensuremath{\mathbf{C}}}
\def\D{\ensuremath{\mathbf{D}}}

\def\AA{\ensuremath{\mathbf{A}^{1,\mathrm{an}}_{A}}}

\def\p{\ensuremath{\mathfrak{p}}}

\def\As{\ensuremath{\mathscr{A}}}
\def\Bs{\ensuremath{\mathscr{B}}}
\def\Cs{\ensuremath{\mathscr{C}}}

\def\Fs{\ensuremath{\mathscr{F}}}
\def\Gs{\ensuremath{\mathscr{G}}}
\def\Hs{\ensuremath{\mathscr{H}}}
\def\Is{\ensuremath{\mathscr{I}}}
\def\Js{\ensuremath{\mathscr{J}}}

\def\Ms{\ensuremath{\mathscr{M}}}

\def\Os{\ensuremath{\mathscr{O}}}

\def\Ss{\ensuremath{\mathscr{S}}}

\def\Us{\ensuremath{\mathscr{U}}}
\def\Vs{\ensuremath{\mathscr{V}}}

\def\Ys{\ensuremath{\mathscr{Y}}}
\def\Zs{\ensuremath{\mathscr{Z}}}

\def\Pc{\ensuremath{\mathcal{P}}}
\def\Qc{\ensuremath{\mathcal{Q}}}
\def\Fc{\ensuremath{\mathcal{F}}}

\def\cn#1#2{\ensuremath{[\![ {#1},{#2}]\!]}}

\def\be{\ensuremath{{\boldsymbol{\varepsilon}}}}

\def\b0{\ensuremath{{\boldsymbol{0}}}}

\begin{document}
\setlength{\baselineskip}{0.55cm}	

\title[Connexit\'e en g\'eom\'etrique analytique $p$-adique]{Un r\' esultat de connexit\' e\\ pour les vari\' et\' es analytiques $p$-adiques.\\ Privil\`ege et noeth\'erianit\'e.}
\author{J\'er\^ome Poineau}
\address{Institut de recherche math\'ematique avanc\'ee, 7, rue Ren\'e Descartes, 67084 Strasbourg, France}
\email{jerome.poineau@math.u-strasbg.fr}

\date{\today}


\begin{abstract}
Soient $k$ un corps ultram\'etrique complet, $X$ un espace $k$-affino\"ide et~$f$ une fonction analytique sur $X$. Notre th\'eor\`eme principal d\'ecrit avec pr\'ecision la variation des composantes connexes g\'eom\'etriques des espaces du type 
$\{x\in X\,|\, |f(x)|\ge\varepsilon\},$
en fonction de $\varepsilon$. Nous obtenons \'egalement des r\'esultats de privil\`ege et de noeth\'erianit\'e au voisinage d'un compact, analogues \`a ceux de la g\'eom\'etrie analytique complexe. 
\end{abstract}

\maketitle

\tableofcontents

\section*{Introduction}\label{Introduction}

Le travail que nous pr\'esentons ici trouve son origine dans une tentative d'adapter au cadre $p$-adique un r\'esultat connu de g\'eom\'etrie analytique complexe. Il s'agit d'un th\'eor\`eme, d\^u \`a J. Frisch, qui assure la noeth\'erianit\'e de l'anneau des germes de fonctions au voisinage de certains compacts. La preuve qui figure dans l'article \cite{Frisch} fait appel, de mani\`ere essentielle, \`a deux notions dont nous ne disposons pas pour des vari\'et\'es analytiques $p$-adiques : celle de voisinage privil\'egi\'e et celle de stratification d'une partie semi-analytique.\\
\indent Cependant, C. B\u anic\u a et O. St\u an\u a\c sil\u a ont abord\'e le probl\`eme de fa\c con l\'eg\`erement diff\'erente et r\'edig\'e, dans \cite{BS2}, 5, fin du \S 3, une d\'emonstration, dont les arguments peuvent s'adapter, sans peine. Nous proposons, en appendice \`a ce texte, un \'enonc\'e du th\'eor\`eme, dans le cadre des espaces d\'efinis sur un corps ultram\'etrique complet, accompagn\'e d'une preuve, calqu\'ee sur la leur.\\ 
\indent La d\'emonstration originale de J. Frisch, bien qu'\`a pr\'esent obsol\`ete, nous a conduit \`a nous int\'eresser aux voisinages privil\'egi\'es. Nous avons alors cherch\'e \`a \'etendre au cadre des vari\'et\'es $p$-adiques le r\'esultat d'A. Douady (\emph{cf.} \cite{Douady}, \S 6, th\'eor\`eme 1) assurant l'existence de voisinages compacts privil\'egi\'es pour les faisceaux coh\'erents. Ainsi que nous l'expliquerons dans la derni\`ere partie de ce texte, nous pouvons en proposer une d\'emonstration fort simple, pour peu que nous disposions d'une sorte de g\'en\'eralisation du th\'eor\`eme d'extension de Riemann. Elle s'\'enonce comme suit :

\begin{thmintro}\label{connexe}
Soit $X$ un espace $k$-affino\"ide ir\-r\'e\-duc\-ti\-ble et $f_{1},\ldots, f_{n}$, avec $n\in\N^*$, des fonctions analytiques sur $X$. Alors il existe un voisinage $V$ de $0$ dans $\R_{+}^n$ tel que le domaine analytique de $X$ d\'efini par 
$$V_{\be}=\bigcup_{1\le j\le n} \{x\in X\,|\, |f_{j}(x)|\ge\varepsilon_{j}\}$$ 
est irr\'eductible, d\`es que le $n$-uplet $\be=(\varepsilon_{1},\ldots,\varepsilon_{n})$ appartient \`a $V$. 
\end{thmintro}

\indent Dor\'enavant, c'est \`a ce probl\`eme que nous consacrerons notre attention. Rappelons que, dans le cadre de la g\'eom\'etrie analytique complexe, le th\'eor\`eme d'extension de Riemann assure, en particulier, qu'un espace irr\'eductible le reste lorsque l'on lui retire un ferm\'e analytique strict. L'analogue de ce th\'eor\`eme pour les vari\'et\'es analytiques $p$-adiques est \'egalement connu de longue date (\emph{cf.} \cite{Bart}, \cite{Lu}).  Dans quelle mesure est-il possible d'\^oter un voisinage d'un tel ferm\'e sans nuire \`a l'irr\'eductibilit\'e?\\
\indent Remarquons que, sur un ouvert, cadre naturel de la g\'eom\'etrie analytique complexe, toute tentative en ce sens serait vou\'ee \`a l'\'echec.
Pour nous en convaincre, consid\'erons l'ouvert du plan complexe d\'efini par 
$$U=\left\{z\in\C\,|\, |\textrm{Im}(z)|<\frac{1}{|\textrm{Re}(z)|+1}\right\}.$$
Alors, d\`es que $\varepsilon$ est assez petit, le domaine d\'efini par $\{z\in U\,|\, |\sin(z)|\ge\varepsilon\}$
poss\`ede une infinit\'e de composantes connexes.\\
\indent N\'eanmoins, le probl\`eme garde un int\'er\^et pour les espaces analytiques d\'efinis sur un corps ultram\'etrique complet. En effet, les mod\`eles locaux de ces derniers, appel\'es espaces affino\"ides, se comportent, \`a bien des \'egards, comme des espaces compacts. Entre autres propri\'et\'es, ils sont quasi-compacts, d\'efinis par un nombre fini d'in\'egalit\'es larges et v\'erifient un principe du maximum.\\
\indent Nous nous placerons donc d\'esormais sur un espace affino\"ide irr\'eductible $X$ d\'efini sur un corps ultram\'etrique complet $k$. Notons $\sqrt{|k^*|}$ le $\Q$-espace vectoriel engendr\'e par le groupe des valeurs du corps de base. Dans le cadre de la g\'eom\'etrie analytique rigide, les fonctions prennent leurs valeurs dans l'ensemble, peu rago\^utant, $\sqrt{|k^*|}\cup\{0\}$ et la topologie des espaces est une topologie de Grothendieck, qui n'est gu\`ere ais\'ee \`a manipuler. L'espace $X$ est quasi-compact, mais, en g\'en\'eral, nous ne pouvons assurer que l'ensemble $V_{\be}$ v\'erifie cette propri\'et\'e que dans le cas o\`u $\be \in (\sqrt{|k^*|})^n$.\\
\indent Il y a de cela une vingtaine d'ann\'ees, V. G. Berkovich proposa une nouvelle approche des vari\'et\'es analytiques sur un corps ultram\'etrique complet (\emph{cf.} \cite{rouge} et \cite{bleu}). La construction qu'il mit en \oe uvre pr\'esentait plusieurs avantages et notamment celui de fournir des espaces poss\'edant de nombreuses propri\'et\'es topologiques remarquables : \`a titre d'exemple, citons la s\'eparation, la compacit\'e et la connexit\'e par arcs locales. Cette derni\`ere propri\'et\'e nous sera fort utile : le r\'esultat que nous avons en vue imposant \`a un certain espace d'\^etre connexe, cela nous simplifiera la t\^ache de pouvoir tracer des chemins.\\
\indent Ajoutons que, dans ce nouveau cadre, les fonctions prennent leurs valeurs dans un ensemble continu, que l'espace $X$ est compact et qu'il en est de m\^eme pour l'espace $V_{\be}$, quel que soit $\be\in\R_{+}^n$. Ces diff\'erentes raisons nous conduisent \`a nous placer, tout au long de ce texte et sans plus le pr\'eciser d\'esormais, dans le cadre des espaces analytiques au sens de V. G. Berkovich. Nous y d\'emontrerons le th\'eor\`eme annonc\'e. Pr\'ecisons, cependant, que le r\'esultat reste valable dans le cadre de la g\'eom\'etrie analytique rigide.\\ 
\indent Restreignons-nous, \`a pr\'esent, au cas d'une seule fonction $f$ sur $X$ et supposons qu'elle soit born\'ee par 1. Consid\'erons les domaines affino\"ides de $X$ de la forme
$$V_{\varepsilon}=\{x\in X\,|\, |f(x)|\ge \varepsilon\},$$
avec $\varepsilon\in [0,1]$. Notre second th\'eor\`eme assure que les composantes connexes des domaines affino\"{\i}des pr\'ec\'edents varient de fa\c{c}on mod\'er\'ee en fonction du param\`etre $\varepsilon$. Pour le d\'emontrer, nous utiliserons, de mani\`ere essentielle, une autre sp\'ecificit\'e des espaces construits par V. G. Berkovich : dans les bons cas, ils se r\'etractent sur un sous-ensemble ferm\'e, appel\'e squelette, qui est muni d'une structure lin\'eaire par morceaux. En particulier, nous parviendrons \`a lire le param\`etre $\varepsilon$ sur un segment r\'eel, hom\'eomorphe \`a $[0,1]$, trac\'e sur le disque unit\'e de dimension $1$.\\
\indent \'Enon\c{c}ons pr\'ecis\'ement le th\'eor\`eme en question. Soit $\bar{k}$ une cl\^oture alg\'ebrique $k$. Nous noterons $\pi_0^g$ le foncteur, d\'efini de la cat\'egorie des espaces $k$-analytiques dans celle des ensembles munis d'une action du groupe d'automorphismes Aut($\bar{k}/k$), qui associe \`a un espace $k$-analytique l'ensemble de ses composantes connexes g\'eom\'etriques.  
\begin{thmintro} \label{partition}
Soient $k$ un corps ultram\'etrique complet, $X$ un espace $k$-affino\"ide et $f$ une fonction analytique sur $X$. Alors il existe une partition finie $\Pc$ de $\R^+$ de la forme
$$\Pc=\{[0,a_0],\mathopen]a_0,a_1],\ldots,\mathopen]a_{r-1},a_r],\mathopen]a_r,+\infty\mathclose[\},$$
o\`u $r\in\N$ et $(a_i)_{0\le i\le r}$ est une suite croissante d'\'el\'ements de $R_{X}\cup\{0\}$, satisfaisant la condition suivante : quel que soit $I\in\Pc$, quels que soient $\varepsilon', \varepsilon \in I$, avec $\varepsilon'\le \varepsilon$, l'inclusion 
$$\{x\in X\,|\, |f(x)|\ge \varepsilon\} \subset \{x\in X\,|\, |f(x)|\ge \varepsilon'\}$$
induit une bijection 
$$\pi_0^{g}(\{x\in X\,|\, |f(x)|\ge \varepsilon\}) \to \pi_0^{g}(\{x\in X\,|\, |f(x)|\ge \varepsilon'\}).$$
Le m\^eme r\'esultat vaut pour le foncteur qui associe \`a un espace $k$-analytique le Aut($\bar{k}/k$)-ensemble de ses composantes irr\'eductibles g\'eom\'etriques.
\end{thmintro}
Dans ce th\'eor\`eme, l'ensemble $R_{X}$ d\'esigne le sous-$\Q$-espace vectoriel de $\R_{+}^*$ engendr\'e par les valeurs non nulles de la norme spectrale sur l'alg\`ebre de $X$. Par exemple, si l'espace $X$ est strictement $k$-affino\"{\i}de, on a $R_{X}=\sqrt{|k^*|}$.

Remarquons que, dans le cadre de la g\'eom\'etrie rigide, A. Abbes et T. Saito (\emph{cf.} \cite{Abbes}, 5.1) ont d\'ej\`a d\'emontr\'e ce dernier r\'esultat pour un intervalle du type $[a,+\infty[$, avec $a>0$, et en ne s'int\'eressant qu'au cardinal de l'ensemble des composantes connexes g\'eom\'etriques. Nous signalons \'egalement que, dans leur article, les bornes des intervalles sont interpr\'et\'es comme les sauts d'une certaine filtration de ramification. \\

La d\'emonstration que nous proposons s'effectuera en quatre \'etapes, correspondant aux quatre premi\`eres parties de ce texte. Dans un premier temps, nous \'etudierons la mani\`ere dont varient les composantes g\'eom\'etriquement connexes des fibres d'un morphisme entre espaces affino\"ides. Lorsqu'elles se r\'ealiseront comme composantes connexes, nous chercherons \`a les rep\'erer par des sections. Des questions proches ont d\'ej\`a \'et\'e trait\'ees pour des morphismes entre sch\'emas : nous savons, par exemple, d'apr\`es \cite{EGAIV3}, 9.7.9, que, pour un morphisme de pr\'esentation finie, la fonction qui a un point de la base associe le nombre g\'eom\'etrique de composantes connexes de sa fibre est localement constructible sur la base. Nous parviendrons \`a nos fins en appliquant des r\'esultats de ce type sur la fibre sp\'eciale d'un mod\`ele formel, judicieusement choisi, du morphisme entre espaces affino\"ides dont nous sommes partis.\\
\indent La difficult\'e principale tient dans la d\'emonstration de l'existence d'un mod\`ele pos\-s\'e\-dant de bonnes propri\'et\'es. Elle nous est assur\'ee par le th\'eor\`eme de la fibre r\'eduite (\emph{cf.} \cite{FIV}), pourvu que le morphisme entre espaces affino\"ides soit plat et \`a fibres g\'eo\-m\'e\-tri\-que\-ment r\'eduites. Dans ce cas, nous parviendrons \`a exhiber une partition finie de la base en domaines analytiques sur lesquels le nombre g\'eom\'etrique de composantes connexes des fibres est constant.\\
\indent Les deuxi\`eme et troisi\`eme parties seront consacr\'ees aux domaines d\'efinis par $$\{x\in X\,|\, |f(x)|\ge\varepsilon\},$$ pour $\varepsilon>0$, dans le cas particulier d'un espace $X$ strictement affino\"ide int\`egre d\'efini sur un corps alg\'ebriquement clos et d'une fonction $f$ de norme spectrale \'egale \`a 1. Nous montrerons que la variation de leurs composantes connexes, en fonction de $\varepsilon$, est li\'ee \`a un probl\`eme du type pr\'ec\'edent. \`A cet effet, nous construirons explicitement un morphisme $\tau$ au-dessus du disque analytique $\D=\Ms(k\{U\})$ de dimension $1$ dont les fibres seront isomorphes, apr\`es extension du corps de base, aux domaines affino\"ides en question. Le param\`etre r\'eel $\varepsilon$ recevra, lui aussi, une interpr\'etation g\'eom\'etrique en termes de valeur absolue de l'\'evaluation de la fonction $U$ sur le disque.\\
\indent Afin d'appliquer les r\'esultats du d\'ebut, nous devrons nous assurer que le morphisme~$\tau$ v\'erifie certaines propri\'et\'es. Nous d\'emontrerons sans peine qu'il est plat, mais buterons sur le caract\`ere g\'eom\'etriquement r\'eduit de l'une des fibres. Dans la troisi\`eme partie, nous modifierons le morphisme $\tau$ de fa\c con \`a passer outre ce probl\`eme. Les techniques mises en jeu rel\`everont, cette fois-ci, de la g\'eom\'etrie alg\'ebrique, puisque nous travaillerons sur des spectres, au sens sch\'ematique, d'alg\`ebres affino\"ides. L'argument principal que nous utiliserons sera le th\'eor\`eme d'\'elimination de la ramification d\'emontr\'e par H. Epp dans \cite{Epp}. Par ce biais, nous parviendrons \`a obtenir des informations sur les fibres de $\tau$ voisines de celle qui pr\'esente des multiplicit\'es et \`a ramener le probl\`eme de la connexit\'e de $\{x\in X\,|\, |f(x)|\ge\varepsilon\}$, pour $\varepsilon$ proche de $0$, \`a celui de $\{x\in X\,|\, |f(x)|>0\}$. Nous concluerons gr\^ace \`a l'analogue ultram\'etrique du th\'eor\`eme de Hartogs (\emph{cf.} \cite{Bart} ou \cite{Lu}).\\
\indent Dans la quatri\`eme partie, nous expliquerons comment d\'eduire les th\'eor\`emes \ref{connexe} et \ref{partition} en toute g\'en\'eralit\'e, \`a partir des cas particuliers consid\'er\'es dans les parties pr\'ec\'edentes.\\

Je tiens \`a remercier Antoine Chambert-Loir pour les nombreux conseils qu'il m'a prodigu\'es. C'est sur ses indications que je me suis int\'eress\'e au th\'eor\`eme de \mbox{H. Epp}, sans lequel ce travail n'aurait pu \^etre men\'e \`a terme. Ma gratitude va \'egalement \`a Antoine Ducros pour avoir suivi avec attention l'avanc\'ee de mes recherches et avoir toujours accept\'e de partager avec moi sa passion pour les espaces de Berkovich. J'exprime \'egalement mes remerciements \`a Ahmed Abbes pour l'int\'er\^et qu'il a port\'e \`a mon travail et ses remarques qui m'ont permis de pr\'eciser le r\'esultat du th\'eor\`eme \ref{partition}. Finalement, je sais gr\'e \`a Qing Liu qui a lu attentivement ce texte et m'a invit\'e \`a \'etendre le r\'esultat du th\'eor\`eme \ref{connexe}.

\section{Connexit\'e des fibres d'un morphisme}

\indent Dans cette partie, nous fixerons un corps ultram\'etrique complet $k$ dont nous supposerons que la valuation n'est pas triviale. Nous noterons $k^\circ$ son anneau de valuation et $\tilde{k}$ son corps r\'esiduel.\\ 
\indent Soit $\Bs$ une alg\`ebre strictement $k$-affino\"ide. Rappelons qu'il existe deux fa\c cons de r\'eduire l'espace affino\"ide $\Ms(\Bs)$ en une vari\' et\' e alg\' ebrique. La premi\`ere, que l'on trouvera expliqu\' ee, par exemple, dans \cite{rouge}, 2.4, utilise la semi-norme spectrale, not\' ee $|.|_{\textrm{sup}}$, sur l'alg\`ebre strictement $k$-affino\"ide $\Bs$. Elle associe \`a l'espace affino\"ide $\Ms(\Bs)$ la vari\'et\'e alg\'ebrique $\textrm{Spec}(\tilde{\Bs})$ o\`u $\tilde{\Bs}$ d\' esigne le quotient de l'anneau $\Bs^{\circ} = \{g\in\Bs\, /\, |g|_{\textrm{sup}}\le 1\}$ par l'id\' eal $\Bs^{\circ\circ} = \{g\in\Bs\, /\, |g|_{\textrm{sup}} < 1\}$. Dans ce cas, l'application de r\'eduction est surjective, anticontinue et induit une bijection entre les composantes connexes.\\
\indent La seconde r\'eduction, due \`a M. Raynaud (\emph{cf.} \cite{tableronde}), consiste \`a interpr\' eter l'espace affino\"ide comme la fibre g\' en\' erique d'un sch\' ema formel plat sur un anneau de valuation. La vari\' et\' e alg\' ebrique associ\'ee est alors d\'efinie comme la fibre sp\' eciale du mod\`ele. On d\'emontre ais\'ement que tout espace affino\"ide admet un mod\`ele formel. Un r\'esultat plus difficile assure m\^eme que tout morphisme peut se r\'ealiser comme un morphisme entre mod\`eles, ce dernier pouvant \^etre choisi plat lorsque le morphisme de d\'epart l'est. Pour plus de d\'etails, nous renvoyons aux articles de r\' ef\' erence \cite{FI} et \cite{FII}.\\

\indent Rappelons qu'un $k^\circ$-sch\'ema formel est dit admissible s'il est localement topologiquement de pr\'esentation finie et plat sur Spf($k^\circ$) (\emph{i.e.} sans $k^\circ$-torsion). Nous dirons qu'un morphisme entre espaces $k$-affino\"ides $\Ms(\mathscr{C})\to\Ms(\Bs)$ est plat lorsque le morphisme associ\'e $\Bs\to\mathscr{C}$ entre alg\`ebres $k$-affino\"ides l'est.

\begin{lem}\label{fplat} 
Tout morphisme plat entre $k^\circ$-sch\'emas formels admissibles le reste apr\`es changement de base par un morphisme entre $k^\circ$-sch\'emas formels admissibles et apr\`es extension des scalaires de $k^\circ$ \`a $L^\circ$, o\`u $L^\circ$ d\'esigne l'anneau de valuation d'une extension ultram\'etrique compl\`ete $L$ de $k$.\\
\indent Ces r\'esultats restent valables pour les morphismes entre espaces strictement $k$-affino\"ides.
\end{lem}
\begin{proof}
Notons $m$ l'id\'eal maximal de $k^\circ$. Soit $\varphi : A\to B$ un morphisme de $k^\circ$-alg\`ebres topologiquement de pr\'esentation finie. D'apr\`es \cite{FI}, 1.6,  le morphisme $\varphi$ est plat si, et seulement si, quel que soit $n\in\N^*$, le morphisme induit $$A\otimes_{k^\circ} (k^\circ/m^n) \to B\otimes_{k^\circ} (k^\circ/m^n)$$
est plat. Le r\'esultat pour les sch\'emas formels en d\'ecoule imm\'ediatement. Celui pour les espaces strictement $k$-affino\"ides s'y ram\`ene puisqu'un morphisme plat entre tels espaces admet un mod\`ele plat, d'apr\`es \cite{FII}, 5.10. 
\end{proof}

Les deux lemmes suivants illustrent l'importance des morphismes plats entre mod\`eles formels. Lorsque, par la suite, nous consid\'ererons un $k^\circ$-sch\'ema formel admissible, nous d\'esignerons sa fibre g\'en\'erique (resp. sp\'eciale) par le m\^eme symbole, auquel nous ajouterons un $\eta$ (resp. une $s$) en indice. Nous adopterons la m\^eme convention pour les morphismes entre tels objets.

\begin{lem}
Sur un sch\'ema formel admissible, l'application de r\'eduction est surjective. 
\end{lem}
\begin{proof} Soit $\Ys$ un $k^\circ$-sch\'ema formel admissible. Nous pouvons le supposer affine, d'alg\`ebre $\Bs$. Soit $\tilde{y}$ un point de la fibre sp\'eciale $\Ys_s$ de $\Ys$. Notons $\tilde{k}(\tilde{y})$ son corps r\'esiduel. Puisque $\Ys_s$ est un $\tilde{k}$-sch\'ema de type fini, il existe une base de transcendance $(T_1,\ldots,T_r)$, avec $r\in\N$, de $\tilde{k}(\tilde{y})$ sur $\tilde{k}$.\\
\indent Soit $K$ le compl\'et\'e du corps $k(U_1,\ldots,U_r)$ pour la norme de {Gau\ss}. Son corps r\'esiduel est isomorphe \`a $\tilde{k}(T_1,\ldots,T_r)$. Consid\'erons le compl\'et\'e $L$ d'une cl\^oture alg\'ebrique de $K$. D'apr\`es \cite{BGR}, 3.4.1/5, son corps r\'esiduel $\tilde{L}$ est une cl\^oture alg\'ebrique de $\tilde{k}(T_1,\ldots,T_r)$ et contient donc un corps isomorphe \`a $\tilde{k}(\tilde{y})$.\\
\indent D'apr\`es le lemme \ref{fplat}, le morphisme 
$$\varphi : \Zs = \Ys\times_{\textrm{Spf}(k^\circ)} \textrm{Spf}(L^\circ) \to \textrm{Spf}(L^\circ)$$ 
est plat, autrement dit, $\Zs$ d\'efinit un $L^\circ$- sch\'ema formel admissible. Par construction, la fibre $\varphi_s^{-1}(\tilde{y})$ au-dessus de $\tilde{y}$ poss\`ede un point ferm\'e de corps r\'esiduel isomorphe \`a $\tilde{L}$. Ce point est encore ferm\'e dans $\Zs_s$. D'apr\`es \cite{Ber}, 1.1.5, c'est l'image d'un point de la fibre g\'en\'erique, ce qui permet de conclure.  
\end{proof}

De ce r\'esultat, d\'eduisons-en un autre, que nous utiliserons \`a de nombreuses reprises :

\begin{lem}\label{plat}
Soit $\varphi : \Ys\to\Zs$ un morphisme de $k^\circ$-sch\'emas formels admissibles. Nous noterons indiff\'eremment $\pi$ les deux morphismes de sp\'ecialisation.
$$\xymatrix{
\Ys_\eta \ar[d]_{\varphi_\eta}\ar^{\pi}[r] &    \ar[d]_{\varphi_s}  \Ys_s\\
\Zs_\eta  \ar[r]^{\pi}& \Zs_s}$$ 
Si le morphisme $\varphi$ est plat, alors les sp\'ecialisations des fibres de $\varphi_\eta$ co\"incident avec les fibres du morphisme sp\'ecialis\'e $\varphi_s$. Autrement dit, quel que soit $z\in\Zs_\eta$, on dispose d'un isomorphisme
$$\pi(\varphi_\eta^{-1}(z))\simeq \varphi_s^{-1}(\pi(z)).$$ 
\end{lem}
\begin{proof}
Seule la surjectivit\'e n\'ecessite une d\'emonstration. Nous pouvons supposer que le sch\'ema formel $\Zs$ est affine, d'alg\`ebre $\Cs$. Soit $z$ un point de $\Zs_\eta$. Il lui correspond un caract\`ere $\chi_z : \Cs\otimes_{k^\circ} k\to \Hs(z)$. On en d\'eduit un morphisme $\chi : \Cs\to \Hs(z)^\circ$, la notation $\Hs(z)^\circ$ d\'esignant l'anneau de valuation du corps ultram\'etrique complet $\Hs(z)$.\\ 
\indent Par le m\^eme raisonnement que dans le lemme pr\'ec\'edent, on d\'emontre que le sch\'ema formel $\Ys\times_{\Zs}\textrm{Spf}(\Hs(z)^\circ)$ est un $\Hs(z)^\circ$-sch\'ema formel admissible. Ses fibres g\'en\'erique et sp\'eciale sont respectivement isomorphes \`a $\varphi_\eta^{-1}(z)$ et $\varphi_s^{-1}(\pi(z))\otimes_{\tilde{k}} \tilde{\Hs(z)}$. Le morphisme $\varphi_s^{-1}(\pi(z))\otimes_{\tilde{k}} \tilde{\Hs(z)} \to \varphi_s^{-1}(\pi(z))$ est surjectif et le lemme pr\'ec\'edent nous permet de conclure.
\end{proof}

Des liens, d\' etaill\' es dans \cite{FIV}, \S 1, existent parfois entre la r\'eduction d\'efinie par la norme spectrale et celle au sens de Raynaud. Citons qu'elles co\"incident lorsque la fibre sp\'eciale du mod\`ele formel est r\'eduite. La prochaine proposition d\'ecoulera de ce r\'esultat. \'Enon\c cons, au pr\'ealable, quelques d\'efinitions.

\begin{defi}
Soit $p : A \to B$ une application continue entre espaces topologiques. Soient un entier $r\in\N$, une partie $P$ de $B$, une famille d'espaces topologiques $(Q_i)_{1\le i\le r}$ et deux familles d'applications continues 
$$\mathbf{s}=(s_i : Q_i \to B)_{1\le i\le r}\quad \textrm{ et }\quad \mathbf{t}=(t_i : Q_i \to A)_{1\le i\le r}.$$ 
Nous dirons que les familles $\mathbf{s}$ et $\mathbf{t}$ \textbf{rep\`erent les composantes connexes des fibres} de $p$ au-dessus de $P$ si les conditions suivantes sont satisfaites :
\begin{enumerate}[a)]
\item quel que soit $i\in\cn{1}{r}$, l'image de l'application $s_i$ recouvre $P$ ;
\item quel que soit $b\in P$, la fibre $p^{-1}(b)$ poss\`ede exactement $r$ composantes connexes ;
\item quels que soient $b\in P$ et la composante connexe $C$ de $p^{-1}(b)$, il existe un unique $i\in\cn{1}{r}$ pour lequel on ait 
$$\forall b'\in s_i^{-1}(b),\, t_i(b')\in C.$$
\end{enumerate}
\end{defi}

\begin{defi}
Soient $\varphi : Y \to Z$ un morphisme de sch\'emas et $P$ une partie de l'espace topologique sous-jacent \`a $Z$. Nous dirons que le morphisme $\varphi$ admet un \textbf{d\'ecoupage} au-dessus de $P$ s'il existe un entier $r\in\N$ et des familles finies 
$$\mathbf{s}=(s_i : Z_i \to Z)_{1\le i\le r}\quad \textrm{ et } \quad \mathbf{t}=(t_i : Z_i\to Z_i\times_Z Y)_{1\le i\le r}$$ de morphismes entre sch\'emas v\'erifiant les conditions suivantes :
\begin{enumerate}[a)]
\item quel que soit $i\in\cn{1}{r}$, le morphisme $t_i$ d\'efinit une section du morphisme $Z_i \times_{Z} Y\to Z_i$, obtenu \`a partir de $\varphi$ par le changement de base $s_i$ ;
$$\xymatrix{
Z_i \times_{Z} Y \ar[d] \ar[r]^{\quad s_i'} & Y \ar[d]_\varphi\\
 Z_i \ar@/^1pc/[u]^{t_i}\ar[r]^{s_i} & Z
 }$$
\item quel que soit $i\in\cn{1}{r}$, le morphisme $s_i$ est \'etale ;
\item les familles $\mathbf{s}$ et $\mathbf{t'}=(s_i' \circ t_i)_{1\le i\le r}$ rep\`erent les composantes connexes des fibres de $\varphi$ au-dessus de $P$. 
\end{enumerate}

\indent Nous adoptons la m\^eme d\'efinition pour un morphisme entre espaces $k$-analytiques en rempla\c cant les morphismes \'etales par des morphismes quasi-\'etales\,\footnote{Les morphismes quasi-\'etales d\'efinis par V. G. Berkovich (\emph{cf.} \cite{VC}, \S 3) correspondent aux morphismes rig-\'etales de la g\'eom\'etrie rigide (\emph{cf.} \cite{FIII}, 3.1).} dont la source est un espace $k$-analytique compact.  
\end{defi}

\begin{defi}
Nous dirons qu'un sch\'ema (resp. espace $k$-analytique) est \textbf{d\'eploy\'e} lorsque ses composantes connexes sont g\'eom\'etriquement connexes.
\end{defi}

\begin{prop}\label{fibreformel}
Soit $\varphi : \Ys \to \Zs$ un morphisme plat entre $k^\circ$-sch\'emas formels admissibles et quasi-compacts. Supposons que le morphisme $\varphi_s$, induit entre les fibres sp\'eciales, soit surjectif et que ses fibres soient g\'eom\'etriquement r\'eduites et d\'eploy\'ees. Alors il existe une partition finie $\Pc$ de $\Zs_s$ v\'erifiant les conditions suivantes :
\begin{enumerate}[a)]
\item les \'el\'ements de $\Pc$ sont des parties constructibles de $\Zs_s$ ;
\item quel que soit $P\in\Pc$, le morphisme $\varphi_\eta$ admet un d\'ecoupage au-dessus du tube de $P$.
\end{enumerate}
\end{prop}
\begin{proof}
Int\'eressons-nous au morphisme $\varphi_s : \Ys_s \to \Zs_s$ entre sch\'emas de type fini sur le corps $\tilde{k}$. Soit $Z$ un ferm\'e irr\'eductible de $\Zs_s$ de point g\'en\'erique $\zeta$. Notons $$\psi : Y = Z\times_{\Zs_s} \Ys_s \to Z$$ le morphisme induit par $\varphi_s$ au-dessus de $Z$. Soient $C_1,\ldots,C_r$, avec $r\in\N$, les composantes connexes de la fibre g\'en\'erique $\psi^{-1}(\zeta)$. Pour chaque $i\in\cn{1}{r}$, choisissons un voisinage ouvert $U_i$ de $C_i$ dans $Y$ dont la trace sur $\psi^{-1}(\zeta)$ soit \'egale \`a $C_i$. L'ensemble des points de $Y$ appartenant \`a au moins deux de ces voisinages forme une partie constructible donc, d'apr\`es le th\'eor\`eme de Chevalley (\emph{cf.} \cite{EGAIV1}, 1.8.4), son image $V$ d\'efinit une partie constructible de $Z$. Puisque la partie $V$ ne contient pas le point g\'en\'erique $\zeta$, elle \'evite m\^eme un ouvert $W$ autour de ce point.\\ 
\indent Une nouvelle utilisation du th\'eor\`eme de Chevalley nous montre qu'il existe un voisinage ouvert $W'$ de $\zeta$ dans $W$ tel que, quel que soit $z\in W'$, les traces des ouverts $U_i$, avec $i\in\cn{1}{r}$, recouvrent la fibre $\psi^{-1}(z)$. Pour $z\in W'$, elles sont donc r\'eunions de composantes connexes disjointes de $\psi^{-1}(z)$. D'apr\`es \cite{EGAIV3}, 9.7.8, le nombre g\'eom\'etrique de composantes connexes des fibres de $\psi$ est constant sur un voisinage ouvert $W''$ de $\zeta$ dans $W'$. Pour $i\in\cn{1}{r}$, notons $V_i$ le voisinage ouvert de $C_i$ d\'efini par $V_i = U_i \cap \psi^{-1}(W'')$. Puisque les fibres de $\varphi_s$ sont d\'eploy\'ees, quel que soit $z\in W''$, les traces des ouverts $V_i$, pour $i\in\cn{1}{r}$, sur la fibre $\psi^{-1}(z)\simeq \varphi_s^{-1}(z)$ sont exactement les composantes connexes de cette fibre.\\ 
\indent Soit $i\in\cn{1}{r}$. Choisissons un ouvert $V'_i$ de $\Ys_s$ dont la trace sur $Y$ soit \'egale \`a $V_i$. Notons $\varphi_i : V'_i \to \Zs_s$ le morphisme induit par $\varphi_s$ sur $V'_i$. Par hypoth\`ese, la fibre $$\varphi_i^{-1}(\zeta)\simeq C_i$$ est g\'eo\-m\'e\-tri\-que\-ment r\'eduite et non vide. Elle contient donc un point $y_i$ en lequel elle est lisse. Puisque le morphisme $\varphi_s$ est plat, ce point est encore lisse dans $\Ys_s$. Quitte \`a restreindre $V'_i$, nous pouvons donc supposer que $\varphi_i$ est lisse. On en d\'eduit qu'il existe un sch\'ema quasi-compact $S_i$, un morphisme \'etale $s_i : S_i\to \Zs_s$ et un point $\zeta'_i$ au-dessus de $\zeta$ tel que le morphisme $S_i \times_{\Zs_s} V'_i\to S_i$, obtenu \`a partir de $\varphi_i$ par le changement de base $s_i$, admette une section $t : S_i\to S_i \times_{\Zs_s} V'_i$, o\`u $t(\zeta'_i)$ s'envoie sur $y_i$.
$$\xymatrix{
S_i \times_{\Zs_s} V'_i \ar[r] \ar[d]& V'_i \ar[d]^{\varphi_i}\\
S_i \ar[r]^s \ar@/^1pc/[u]^{t}& \Zs_s
}$$ 
\indent Les images des morphismes \'etales $s_i$, pour $i\in\cn{1}{r}$, contiennent un voisinage ouvert commun de $\zeta$ dans $Z$. Par construction, le morphisme $\varphi_s$ y admet un d\'ecoupage. Un argument de r\'ecurrence noeth\'erienne nous montre ensuite qu'il existe une partition $\Pc$ de $\Zs_s$ en parties constructibles au-dessus desquelles le morphisme $\varphi_s$ admet un d\'ecoupage.\\
\indent Remarquons que nous pouvons relever les constructions pr\'ec\'edentes aux sch\'emas formels. Consid\'erons, en effet, une restriction $\psi$ de $\varphi_s$ \`a un ouvert lisse $U$ et un morphisme \'etale $s : S \to \Zs_s$ tel que le morphisme $S\times_{\Zs_s} U\to S$, obtenu \`a partir de $\psi$ par le changement de base $s$, admette une section $t$. Notons $\Us$ le sous-sch\'ema formel ouvert de $\Ys$ dont l'espace topologique sous-jacent est l'ouvert $U$. D'apr\`es \cite{VC}, 2.1, le morphisme $s$ admet un mod\`ele formel \'etale $\Ss \to\Zs$. D'autre part, puisque le morphisme $\Ss\times_{\Zs} \Us \to \Ss$ est lisse, la propri\'et\'e de rel\`evement infinit\'esimal nous assure que la section $t$ se rel\`eve en une section formelle $T$.
$$\xymatrix{
\Ss \times_{\Zs} \Us \ar[r] \ar[d]& \Us \ar[d]^{\varphi}\\
\Ss \ar[r]\ar@/^1pc/[u]^{T}& \Zs
}$$ 
\indent Expliquons, \`a pr\'esent, comment passer des sch\'emas formels \`a leur fibre g\'en\'erique. Soient $z$ un point de $\Zs_\eta$ et $\tilde{z}$ son image dans $\Zs_s$ par l'application de sp\'ecialisation. Puisque le morphisme $\varphi$ est plat et \`a fibres g\'eom\'etriquement r\'eduites, la r\'eduction de $\varphi_\eta^{-1}(z)$ obtenue par la norme spectrale co\"incide avec celle au sens des mod\`eles, \`a savoir la fibre $\varphi_s^{-1}(\tilde{z})$. D'apr\`es les propri\'et\'es de la r\'eduction par la norme spectrale, chaque composante connexe de la fibre analytique est le tube d'une composante connexe de sa r\'eduction, et \emph{vice versa}.\\
\indent D'autre part, si $\Ss \to\Zs$ est un morphisme \'etale entre sch\'emas formels quasi-compacts, alors le morphisme $\Ss_\eta \to\Zs_\eta$, induit entre les fibres g\'en\'eriques, est un morphisme quasi-\'etale entre espaces $k$-analytiques compacts. 
\end{proof}

Le th\'eor\`eme suivant concerne le comportement des composantes connexes des fibres d'un morphisme entre espaces analytiques. La d\'emonstration repose sur le th\'eor\`eme de la fibre r\'eduite (\emph{cf.} \cite{FIV}) qui nous permet de nous ramener \`a un mod\`ele du morphisme satisfaisant les hypoth\`eses de la proposition pr\'ec\'edente. Cette fois encore, nous commen\c cons par poser une d\'efinition.

\begin{defi}
Soit $Y$ un espace $k$-affino\"ide. Une partie $V$ de $Y$ est dite \textbf{simple} si elle peut s'obtenir par combinaison bool\'eenne finie de domaines affino\"ides du type 
$$\{z\in Y\,|\, |h(z)|=1\}$$
o\`u $h$ d\'esigne une fonction analytique sur $Y$ de norme spectrale \'egale \`a $1$.
\end{defi}

Remarquons, d\`es \`a pr\'esent, qu'une partie simple d'un espace $k$-affino\"ide est voisinage de chacun de ses points rigides. En effet, elle s'obtient par r\'eunion et intersection d'un nombre fini de parties qui sont soit des domaines affino\"ides, soit des ouverts et pour lesquelles ce r\'esultat est vrai.    

\begin{thm} \label{fibre}
Soit $\varphi : Y\to Z$ un morphisme plat et surjectif entre espaces strictement $k$-affino\"ides dont les fibres soient g\'eom\'etriquement r\' eduites et d\'eploy\'ees. Alors il existe une partition finie $\Pc$ de $Z$ v\'erifiant les conditions suivantes :
\begin{enumerate}[a)]
\item les \'el\'ements de $\Pc$ sont des parties simples de domaines affino\"ides de $Z$ ;
\item quel que soit $P\in\Pc$, le morphisme $\varphi$ admet un d\'ecoupage au-dessus de $P$.
\end{enumerate}
\end{thm}
\begin{proof}
\indent Les hypoth\`eses de l'\' enonc\' e nous permettent d'appliquer le th\' eor\`eme de la fibre r\' eduite. Celui-ci nous assure qu'il existe un diagramme commutatif de sch\' emas formels 
$$\xymatrix{
& \Ys'\ar_{\chi}[d] \ar@/^3pc/[dd]^{\psi'}\\
\Ys \ar[d]_{\psi} & \Ys\times_{\Zs}\Zs' \ar[l] \ar[d]\\
\Zs & \Zs' \ar[l]}$$ 
o\`u 
\begin{enumerate}[a)]
\item le morphisme $\psi :\Ys\to\Zs$ induit le morphisme $\varphi :Y \to Z$ par passage aux fibres g\'en\'eriques ;
\item le morphisme $Z'\to Z$, o\`u $Z'$ d\' esigne la fibre g\' en\' erique de $\Zs'$, est quasi-\' etale et surjectif ;
\item le morphisme $\chi : \Ys'\to \Ys\times_{\Zs} \Zs'$ est fini et induit un isomorphisme entre les fibres g\'en\'eriques ;
\item le morphisme $\psi' : \Ys'\to \Zs'$ est plat et ses fibres sont g\' eom\' etriquement r\' eduites.
\end{enumerate}

\indent V\'erifions que nous pouvons appliquer la proposition pr\'ec\'edente au morphisme $\psi'$. Les $k^\circ$-sch\'emas formels obtenus comme mod\`eles des espaces analytiques sont bien admissibles et quasi-compacts, en vertu de l'\'equivalence de cat\'egories \cite{FI}, 4.1. Tous les sch\'emas formels que nous consid\'ererons sont \'egalement de ce type.\\
\indent Les fibres de $\psi'_\eta$ sont isomorphes, apr\`es extension du corps de base, \`a des fibres de $\psi_\eta = \varphi$. Par hypoth\`ese, le morphisme $\varphi$ est surjectif. Il en est donc de m\^eme pour $\psi'_\eta$, puis pour $\psi'_s$, d'apr\`es \ref{plat}.\\
\indent Il nous reste \`a d\'emontrer que les fibres de $\psi'$ sont d\'eploy\'ees. Puisque le morphisme $\chi$ est fini, il suffit m\^eme de le v\'erifier sur le morphisme $\Ys\times_{\Zs} \Zs' \to \Zs'$, autrement dit sur le morphisme $\psi$. Or le morphisme $\psi$ est plat et \`a fibres g\'eom\'etriquement r\'eduites, donc le nombre de composantes connexes (resp. g\'eom\'etriquement connexes) des fibres de $\psi_\eta$ est identique \`a celui de leur r\'eduction. D'apr\`es \ref{plat}, toutes les fibres de $\psi_s$ peuvent \^etre obtenues par de telles r\'eductions. Elles sont d\'eploy\'ees, puisque les fibres de $\varphi$ le sont, par hypoth\`ese.\\
\indent D'apr\`es la proposition \ref{fibreformel}, il existe une partition finie $\Pc$ de $\Zs'_s$ en parties constructibles au-dessus des tubes desquelles le morphisme $\psi'_\eta$ admet un d\'ecoupage. Notons $\Qc$ l'ensemble de ces tubes. Puisque le morphisme $\chi$ induit un isomorphisme entre les fibres g\'en\'eriques et que le morphisme $\lambda : Z'\to Z$ est quasi-\'etale, le morphisme $\varphi=\psi_\eta$ admet un d\'ecoupage au-dessus de l'image par $\lambda$ de toute partie de $\Qc$. Puisque le morphisme $Z'\to Z$ est surjectif, nous d\'emontrons ainsi l'existence d'un recouvrement fini de $Z$ par des parties au-dessus desquelles le morphisme $\varphi$ admet un d\'ecoupage. Il est ais\'e d'en d\'eduire une partition de $Z$ v\'erifiant la m\^eme propri\'et\'e et compos\'ee uniquement de combinaisons bool\'eennes des parties pr\'ec\'edentes.\\ 
\indent Pour clore la d\'emonstration, il nous reste \`a v\'erifier que l'image par $\lambda$ de tout \'el\'ement de $\Qc$ est de la forme d\'esir\'ee. Soit $Q\in \Zs'_\eta$ un \'el\'ement de $\Qc$. Il est obtenu comme le tube d'une partie constructible $P$ de $\Zs'_s$. Le morphisme $\lambda : Z'\to Z$ est plat donc, d'apr\`es \cite{FII}, 5.2, il existe un diagramme commutatif de sch\'emas formels 
$$\xymatrix{
\Zs' \ar[r] & \Zs \\
\Zs'_0 \ar^{\psi_0}[u]\ar^{\mu}[r] & \Zs_0  \ar[u]}$$ 
o\`u le morphisme $\mu : \Zs'_0 \to \Zs_0$ est plat et induit encore le morphisme $\lambda : Z'\to Z$ entre les fibres g\'en\'eriques.\\
\indent Le tube $Q\subset  Z'$ de $P$ est identique \`a celui de $(\psi_0)_s^{-1}(P)$.
Le morphisme $\mu_s : (\Zs'_0)_s\to (\Zs_0)_s$ induit par $\mu$ est de type fini et envoie donc la partie constructible $(\psi_0)_s^{-1}(P)$ sur une partie constructible $P_0$ de $(\Zs_0)_s$, d'apr\`es le th\'eor\`eme de Chevalley. Puisque le morphisme $\mu$ est plat, d'apr\`es \ref{plat}, le tube de $P_0$ dans $Z$ n'est autre que l'image de $Q$ dans $Z$ par $\lambda$.\\         
\indent Il nous reste, d\'esormais, \`a montrer que le tube d'une partie constructible de $\Zs_s$ est une r\'eunion finie de parties simples de domaines affino\"ides de $\Zs_\eta$. Ce r\'esultat provient directement des propri\'et\'es de la r\'eduction lorsque le sch\'ema formel consid\'er\'e est affine. Nous concluons en remarquant que le sch\'ema formel quasi-compact $\Zs$ admet un recouvrement fini par de tels ouverts.  
\end{proof}

Rappelons, maintenant, que le nombre g\' eom\' etrique de composantes connexes d'un espace affino\"ide ne change pas lorsque l'on \' etend le corps de base. On trouvera la d\' emonstration de ce r\' esultat dans \cite{D}, th\'eor\`eme 5.5. On pourrait \' egalement le d\' eduire des th\'eor\`emes d'extension du corps de base pour les faisceaux \' etales qui figurent dans \cite{bleu}.\\
\indent Par un raisonnement tr\`es proche de celui que nous venons de mettre en \oe uvre, nous obtenons le th\'eor\`eme suivant :

\begin{thm}\label{variation}
Soit $\varphi : Y\to Z$ un morphisme plat entre espaces strictement $k$-affino\"ides dont les fibres soient g\'eom\'etriquement r\' eduites. Alors il existe une partition finie de $Z$ en parties simples de domaines affino\"ides au-dessus desquelles le nombre g\'eom\'etrique de composantes connexes des fibres est constant. En particulier, ce nombre est constant au voisinage des points rigides.
\end{thm}
\begin{proof}
Appliquons le th\'eor\`eme de la fibre r\'eduite et reprenons les notations de la preuve pr\'ec\'edente. Le morphisme $\psi'_s : \Ys'_s \to \Zs'_s$ est un morphisme de type fini entre deux vari\'et\'es alg\'ebriques de type fini sur $\tilde{k}$. D'apr\`es \cite{EGAIV3}, 9.7.9, il existe donc une partition finie $\Pc$ de $\Zs'_s$ en parties constructibles au-dessus desquelles le nombre g\'eom\'etrique de composantes connexes des fibres de $\psi'_s$ soit constant.\\
\indent Soit $P\in\Pc$. Notons $Q\subset Z'$ le tube de $P$. Puisque le morphisme $\psi'$ est plat et que ses fibres sont g\'eom\'etriquement r\'eduites, le nombre g\'eom\'etrique de composantes connexes des fibres de $\psi'_\eta$ est constant au-dessus de tout point de $Q$. On en d\'eduit le m\^eme r\'esultat pour le morphisme $\psi$ au-dessus de tout point de $\lambda(Q)$. Par un raisonnement en tout point identique \`a celui expos\'e dans la preuve pr\'ec\'edente, nous obtenons une partition de $Z$ en parties simples de domaines affino\"ides jouissant des m\^emes propri\'et\'es.
\end{proof}

\section{Loin de l'hypersurface}

Consacrons-nous, tout d'abord, \`a la d\'emonstration du th\'eor\`eme 2, dans un cas particulier. Dans toute cette partie, $k$ d\'esignera un corps alg\'ebriquement clos dont la valuation n'est pas triviale, $X$ un espace strictement $k$-affino\"ide int\`egre d'alg\`ebre $\As$ et $f$ une fonction analytique sur $X$ de norme spectrale \'egale \`a $1$. Pour $\varepsilon >0$, d\' efinissons le domaine affino\"ide $V_\varepsilon$ de $X$ par 
$$V_\varepsilon = \{x\in X\, /\, |f(x)|\ge \varepsilon\}.$$ 
Nous ne nous int\'eresserons, pour l'instant, qu'aux composantes connexes des espaces~$V_\varepsilon$, avec $\varepsilon>0$.

T\^achons, tout d'abord, de remplacer le param\`etre r\' eel $\varepsilon$ par un autre que nous saurons interpr\' eter g\' eom\' etriquement. Pour ce faire, pla\c cons-nous au-dessus du disque de dimension $1$ et de rayon $1$ d\'efini par $\D=\Ms(k\{U\})$. Dans la suite, nous noterons simplement $0$ le point rigide de $\D$ d\'efini par l'\'equation $U=0$. L'injection $k\{U\}\hookrightarrow \As\{T,U\}$ induit un morphisme $t : k\{U\} \to \As\{T,U\}/(fT-U)$. Nous noterons $\tau$ le morphisme correspondant entre espaces affino\"ides. L'alg\`ebre $k$-affino\"ide \mbox{$\As\{T,U\}/(fT-U)$} h\'erite de nombreuses propri\'et\'es de l'alg\`ebre $k$-affino\"ide $\As$. Le lemme suivant en fournit un exemple.
\begin{lem}\label{iso}
Les alg\`ebres $k$-affino\"ides $\As\{T,U\}/(fT-U)$ et $\As\{T\}$ sont isomorphes. En particulier, l'alg\`ebre $\As\{T,U\}/(fT-U)$ est int\`egre. En outre, elle est int\'egralement close lorsque $\As$ l'est.
\end{lem}
\begin{proof}
On v\'erifie sans peine que le morphisme de $\As$-alg\`ebres  
$$s : \As\{T,U\}/(fT-U)\to \As\{T\}$$ d\'efini par $s(T)=T$ et $s(U)=fT$
est un isomorphisme dont l'inverse est le morphisme de $\As$-alg\`ebres 
$$s^{-1} : \As\{T\}\to \As\{T,U\}/(fT-U)$$
d\'efini par $s^{-1}(T)=T$.\\
\indent Supposons, \`a pr\'esent, que $\As$ soit int\'egralement close. 
Puisque le sch\'ema Spec($\As$) est normal, il en est de m\^eme du sch\'ema Spec($\As[T]$), ainsi que de son analytifi\'e $\mathbf{A}^{1,\textrm{an}}_{\As}$, d'apr\`es \cite{bleu}, 2.2.7. Le domaine affino\"ide $\Ms(\As\{T\})$ de ce dernier est donc, lui aussi, normal, d'apr\`es \cite{bleu}, 2.2.1. Puisque l'anneau $\As\{T\}$ est int\`egre et normal, il est finalement int\'egralement clos.
\end{proof}

Le lemme suivant met en lumi\`ere le lien g\'eom\'etrique recherch\'e :
\begin{lem}
Quel que soit $x\in \D\setminus\{0\}$, la projection 
$$\Ms(\As\{T,U\}/(fT-U)) \to \Ms(\As) = X$$
induit un isomorphisme 
$$\tau^{-1}(x) \xrightarrow{\sim} V_{\varepsilon}\, \hat{\otimes}_k\,  \Hs(x),$$
o\`u $\varepsilon=|U(x)|>0$.\\
\indent En particulier, les fibres du morphisme $\tau$ au-dessus des points de $\D\setminus\{0\}$ sont g\'eom\'etriquement r\'eduites. 
\end{lem}
\begin{proof}
Soit $x\in \D\setminus\{0\}$. Notons $\Bs=\As\hat{\otimes}_k \Hs(x)$. L'alg\`ebre de la fibre $\tau^{-1}(x)$ n'est autre que 
$$(\As\{T,U\}/(fT-U))\hat{\otimes}_{k\{U\}}\, \Hs(x) \simeq \Bs\{T\}/(fT-U(x))\simeq \Bs\{\varepsilon\, T\}/(fT-1),$$
o\`u $\varepsilon=|U(x)|>0$. On reconna\^it l'alg\`ebre du domaine affino\"ide de $\Ms(\Bs)$ d\'efini par 
$$\{y\in\Ms(\Bs)\,|\, |f(y)|\ge\varepsilon\}$$
et qui est isomorphe \`a $V_{\varepsilon}\, \hat{\otimes}_k\,  \Hs(x)$.\\
\indent Passons \`a la seconde partie du lemme. L'espace affino\"ide $X$ est r\'eduit, donc son domaine affino\"ide $V_\varepsilon$ l'est aussi, d'apr\`es \cite{bleu}, 2.2.1. Puisque le corps $k$ est alg\'ebriquement clos, l'espace $V_\varepsilon$ est g\'eom\'etriquement r\'eduit et il en est de m\^eme pour $\tau^{-1}(x) \simeq V_{\varepsilon}\, \hat{\otimes}_k\,  \Hs(x)$, d'apr\`es \cite{D}, 4.17.
\end{proof}

Pour $\varepsilon\in [0,1]$, notons $\eta_{\varepsilon}$ le point de $\D$ associ\' e \`a la valeur absolue d\' efinie par $\sum_{i\ge 0} a_i\, U^i \in k\{U\} \mapsto \max_{i\ge 0}\{|a_i|\, \varepsilon^i\} \in\R_+$. D'apr\`es le lemme, la fibre de $\tau$ au-dessus du point $\eta_{\varepsilon}$ est isomorphe \`a $V_{\varepsilon}\, \hat{\otimes}\,  \Hs(\eta_{\varepsilon})$, quel que soit $\varepsilon\in\, ]0,1]$.\\

Afin de pouvoir appliquer les r\'esultats du paragraphe pr\'ec\'edent, nous avons besoin d'une propri\'et\'e de platitude, que nous d\'emontrons ici. 
\begin{lem}
Le morphisme 
$$\tau : \Ms\left(\As\{T,U\}/(fT-U)\right) \to \Ms(k\{U\})=\D$$ 
est plat.    
\end{lem}
\begin{proof}
L'anneau $k\{U\}$ \'etant principal, il nous suffit de montrer que l'alg\`ebre $\As\{T,U\}/(fT-U)$ ne poss\`ede aucun \' el\' ement de $k\{U\}$-torsion. Il nous suffit m\^eme de montrer que le morphisme
$$t : k\{U\} \to \As\{T,U\}/(fT-U)\simeq \As\{T\}$$ 
est injectif, puisque l'alg\`ebre $\As\{T,U\}/(fT-U)\simeq \As\{T\}$ est int\`egre. L'interpr\'etation g\'eom\'etrique des fibres du morphisme $\tau$ nous montre que sa fibre au point $\eta_1$ n'est pas vide. En particulier, toute fonction $g$ de $k\{U\}$ v\'erifiant $t(g)=0$ est nulle en $\eta_1$ et donc nulle sur $\D$. 
\end{proof}

\begin{rem}
La fibre du morphisme $\tau$ au-dessus du point $0$ de $\D$ est isomorphe \`a l'espace $\Ms(\As\{T\}/(fT))$ et ne saurait donc \^etre r\' eduite lorsque $f$ poss\`ede des multiplicit\' es. Ce probl\`eme fera l'objet du prochain paragraphe.\\
\indent Remarquons n\'eanmoins que si l'hypersurface de $X$ d\'efinie par l'\'equation $f=0$ est r\'eduite, un calcul simple montre que la fibre $\tau^{-1}(0)$ l'est aussi. En outre, elle est connexe, puisqu'elle est r\'eunion de l'espace $X$ et d'une droite au-dessus du lieu d'annulation de $f$ dans $X$. Le th\'eor\`eme \ref{variation} appliqu\'e au morphisme $\tau$ et au point rigide $0$ de $\D$ entra\^ine alors que le domaine affino\"ide $\{x\in X\,|\, |f(x)|\ge \varepsilon\}$ est connexe, d\`es que $\varepsilon$ est assez petit. S'en d\'eduit, en particulier, l'analogue du th\'eor\`eme de Hartogs. 
\end{rem}

Int\'eressons-nous, \`a pr\'esent, \`a la variation des composantes connexes des domaines $V_\varepsilon$, pour $\varepsilon>0$. \'Enon\c cons tout d'abord un lemme. Le caract\`ere fini mis \`a part, nous red\'emontrons ici, dans un cas \'el\'ementaire, le r\'esultat \cite{smoothII}, 6.3.1. Rappelons que nous avons suppos\'e le corps $k$ alg\'ebriquement clos et de valuation non triviale. Par cons\'equent, l'\'egalit\'e $\sqrt{|k^*|}=|k^*|$ est v\'erifi\'ee.

\begin{lem}\label{trace}
La trace d'une partie simple de $\D$ sur le segment
$$\{\eta_\varepsilon,\ 0\le \varepsilon\le 1\}\simeq [0,1]$$
est une r\'eunion finie d'intervalles dont les bornes sont des \'el\'ements de $|k^*|\cup\{0,+\infty\}$. 
\end{lem}
\begin{proof}
Puisqu'une partie simple est obtenue, par d\'efinition, comme une combinaison bool\'eenne finie de domaines affino\"ides, il suffit de d\'emontrer le r\'esultat pour ces derniers. D'apr\`es le th\'eor\`eme de Gerritzen et Grauert (\emph{cf.} \cite{salg}, 2.4), tout domaine affino\"ide de $\D$ peut s'\'ecrire comme r\'eunion finie de domaines rationnels, eux-m\^emes intersections de domaines du type $\{z\in \D\,|\, |g(z)|\le |h(z)|\}$, o\`u $g$ et $h$ d\'esignent des fonctions analytiques sur $\D$. Remarquons encore que, si $g = \sum_{i\in\N} b_i\, U^i$ et $h = \sum_{i\in\N} c_i\, U^i$ dans $k\{U\}$, les fonctions, d\'efinies de $[0,1]$ dans $\R^+$,
$$\varepsilon \mapsto \max_{i\in\N}\, \{\ln(|b_i|)+i\ln(\varepsilon)\}  \quad \textrm{ et } \quad \varepsilon \mapsto \max_{i\in\N}\, \{\ln(|c_i|)+i\ln(\varepsilon)\} $$
sont lin\'eaires par morceaux. Nous concluons gr\^ace \`a l'\'equivalence
$$|g(\eta_\varepsilon)|\le |h(\eta_\varepsilon)|  \Longleftrightarrow \max_{i\in\N}\, \{\ln(|b_i|)+i\ln(\varepsilon)\} \le \max_{i\in\N}\, \{\ln(|c_i|)+i\ln(\varepsilon)\},$$
qui est v\'erifi\'ee quel que soit $\varepsilon\in [0,1]$.\\
\end{proof}

Venons-en au r\'esultat concernant la variation des composantes connexes. Signalons que si l'on ne s'int\'eresse qu'\`a leur nombre, on retrouve un th\'eor\`eme d'A.~Abbes et T.~Saito (\emph{cf.} \cite{Abbes}, 5.1).

\begin{thm}\label{min}
Soit $k$ un corps ultram\'etrique complet alg\'ebriquement clos et dont la valuation n'est pas triviale. Soit $X$ un espace strictement $k$-affino\"ide int\`egre et $f$ une fonction analytique sur $X$ dont la norme spectrale vaut $1$. Soit $m\in\,]0,1]\cap \sqrt{|k^*|}$. Alors il existe une partition finie $\Pc$ de $[m,1]$ en intervalles v\'erifiant la condition suivante : quel que soit $I\in\Pc$, quels que soient $\varepsilon', \varepsilon \in I$, avec $\varepsilon'\le \varepsilon$, l'inclusion 
$$\{x\in X\,|\, |f(x)|\ge \varepsilon\} \subset \{x\in X\,|\, |f(x)|\ge \varepsilon'\}$$
induit une bijection 
$$\pi_0(\{x\in X\,|\, |f(x)|\ge \varepsilon\}) \to \pi_0(\{x\in X\,|\, |f(x)|\ge \varepsilon'\}).$$
En outre, les bornes des intervalles sont des \'el\'ements de $|k^*|\cup\{0,+\infty\}$.
\end{thm}
\begin{proof}
Notons $V$ le domaine strictement $k$-affino\"ide de $\D$ d\'efini par 
$$V = \{z\in\D\,|\, |U(z)|\ge m\}.$$ 
D'apr\`es \ref{fplat}, le morphisme $\tau'$ d\'eduit de $\tau$ par le changement de base $V\hookrightarrow \D$ est plat. Ses fibres, isomorphes, apr\`es extension du corps de base, \`a des espaces du type $V_\varepsilon$, avec $\varepsilon>0$, sont g\'eom\'etriquement r\'eduites et d\'eploy\'ees. D'apr\`es \ref{fibre}, il existe donc une partition finie $\Pc$ de $V$ en parties simples de domaines affino\"ides au-dessus desquelles le morphisme $\tau'$ admet un d\'ecoupage. Consid\'erons l'un des morphismes quasi-\'etales $\varphi : Z \to V$, o\`u $Z$ est un espace $k$-analytique compact, intervenant dans le d\'ecoupage.\\
\indent Le raisonnement qui suit fait intervenir, dans un cas simple, la notion de squelette. On la trouvera introduite dans \cite{smoothI}, \S 5. L'espace $k$-affino\"ide $V$ est isomorphe \`a la fibre g\'en\'erique du $k^\circ$-sch\'ema formel pluristable non d\'eg\'en\'er\'e 
$$\Vs = \textrm{Spf}(k^\circ\{U,V\}/(UV-\alpha)),$$
o\`u $k^\circ$ d\'esigne l'anneau de valuation de $k$ et $\alpha$ un \'el\'ement de $k^\circ$ de valeur absolue $m$. Le squelette $\textrm{S}(\Vs)$ du sch\'ema formel $\Vs$ est le segment $$J = \{\eta_\varepsilon,\ m\le \varepsilon\le 1\}\simeq [m,1],$$
trac\'e sur la fibre g\'en\'erique $\Vs_\eta\simeq V$. D'apr\`es \cite{squelette}, 3.1, il existe alors une unique structure $\sqrt{|k^*|}$-lin\'eaire par morceaux (au sens de \cite{smoothII}, \S 1) sur $\Delta = \varphi^{-1}(\textrm{S}(\Vs))$ telle que l'application
$$\varphi_{|\Delta} : \Delta \to \textrm{S}(\Vs)$$
soit lin\'eaire par morceaux et soit G-localement une immersion. En particulier, puisque $\Delta$ est compact, il existe une partition finie de $\Delta$ en parties lin\'eaires qui sont hom\'eomorphes \`a leur image par $\varphi_{|\Delta}$, elle-m\^eme lin\'eaire. Pour chaque image $Q$, nous pouvons construire, \`a partir de la section associ\'ee \`a $\varphi$, une section de $\tau'$ au-dessus de $Q$.\\
\indent En proc\'edant de m\^eme pour chaque morphisme \'etale, nous obtenons finalement, pour chaque \'el\'ement $P$ de $\Pc$, une partition $\Qc_P$ de $P\cap J$ en un nombre fini de parties lin\'eaires et, au-dessus de chaque $Q\in\Qc_P$, un ensemble fini $\mathcal{T}$ de sections de $\tau'$ au-dessus de $Q$ v\'erifiant la condition suivante : quel que soit $z\in Q$, chaque composante connexe de la fibre $\tau'^{-1}(z)$ contient un et un seul \'el\'ement de la forme $t(z)$, avec $t\in\mathcal{T}$.\\ 
\indent Examinons, \`a pr\'esent, la forme des parties $Q$ consid\'er\'ees pr\'ec\'edemment. Nous souhaitons montrer qu'elles sont r\'eunions finies d'intervalles \`a coordonn\'ees dans $\sqrt{|k^*|}\cup\{0,+\infty\}$. C'est le cas pour les traces des \'el\'ements de $\Pc$ sur $J$, d'apr\`es le lemme \ref{trace}, et donc pour leurs parties $\sqrt{|k^*|}$-lin\'eaires par morceaux.\\
\indent Finalement, les sections de $\tau'$ d\'efinies pr\'ec\'edemment sont d\'efinies sur des intervalles contenus dans $J$. Soient $I$ un tel intervalle et $\varepsilon \in I$. Rappelons que l'image de la fibre de $\tau'$ au-dessus de $\eta_\varepsilon$ par le morphisme 
$$\pi : \Ms(\As\{T,U\}/(fT-U)) \to \Ms(\As) = X $$ 
est isomorphe au domaine affino\"ide $V_\varepsilon$. Pour $\varepsilon'\in I$, $\varepsilon'\le \varepsilon$, les images des sections par le morphisme $\pi$ joignent les composantes connexes de $V_\varepsilon$ \`a celles de $V_\varepsilon'$. Puisque les diff\'erentes sections aboutissent \`a des composantes connexes distinctes et que toutes sont atteintes, on en d\'eduit que les composantes connexes de $V_\varepsilon$ sont les traces de celles de $V_\varepsilon'$.     
\end{proof}

\section{\'Elimination des multiplicit\'es}

Conservons les hypoth\`eses de la partie pr\'ec\'edente : le corps $k$ est un corps al\-g\'e\-bri\-que\-ment clos dont la valuation n'est pas triviale, l'espace $X$ un espace strictement $k$-affino\"ide int\`egre et la norme spectrale de la fonction $f$ vaut $1$. Nous d\'emontrons ici le th\'eor\`eme \ref{connexe} dans ce cas particulier.\\
\indent Remarquons que, puisque le morphisme de normalisation est continu et surjectif, nous pouvons, quitte \`a remplacer $X$ par son normalis\'e, supposer que l'espace $X$ est normal. D'apr\`es \cite{D}, 4.18, les domaines $V_\varepsilon$, pour $\varepsilon>0$, sont alors normaux. Il nous suffit donc de montrer qu'ils sont connexes, pour $\varepsilon$ assez petit.\\ 
\indent D'apr\`es le paragraphe pr\'ec\'edent, nous disposons d'un morphisme plat, 
$$\tau : \Ms(\As\{T,U\}/(fT-U))\to \D,$$ 
dont seule la fibre au-dessus du point $0$ peut pr\'esenter des multiplicit\'es. Nous allons montrer qu'il est possible de modifier ce morphisme de fa\c con que toutes ses fibres deviennent r\'eduites.\\ 
\indent Dans les raisonnements qui suivent, nous quittons le domaine des espaces analytiques pour celui des sch\'emas. Il nous faut donc introduire de nouveaux objets. Soient $D=\textrm{Spec}(k\{U\})$, $F=\textrm{Spec}(\As\{T,U\}/(fT-U))$ et $\alpha : F\to D$ le morphisme induit par $t : k\{U\}\to \As\{T,U\}/(fT-U)$. Soient $x$ un point ferm\'e de $D\setminus\{0\}$ et $\mathbf{x}$ le point rigide de $\D$ qui lui correspond. La fibre de $\alpha$ au-dessus de $x$ a m\^eme anneau que la fibre de $\tau$ au-dessus de $\mathbf{x}$. En particulier, elle est r\'eduite.\\ 
\indent D'apr\`es \ref{iso}, le sch\'ema $F$ est normal. Le probl\`eme de r\'eduction des fibres auquel nous sommes confront\'es se ram\`ene donc \`a un probl\`eme de multiplicit\'es g\'en\'eriques. En effet, une hypersurface principale d'un sch\'ema normal et noeth\'erien est r\'eduite si, et seulement si, elle est g\'en\'eriquement r\'eduite (c'est une cons\'equence de la condition $(\textrm{S}_2)$, \textit{cf.} \cite{M}, 17.I).\\
\indent Rappelons que si $R$ et $S$ sont deux anneaux de valuation discr\`ete et que $S$ domine $R$, on dit que $S$ est faiblement non ramifi\'e au-dessus de $R$ lorsque l'id\'eal maximal de $R$ engendre l'id\'eal maximal de $S$. Si $x$ d\'esigne un point ferm\'e de $D$, la remarque pr\'ec\'edente entra\^ine que la fibre $\alpha^{-1}(x)$ au-dessus de $x$ est r\'eduite si, et seulement si, l'anneau $\Os_{F,\eta}$ est faiblement non ramifi\'e au-dessus de l'anneau $\Os_{D,x}$, pour tout point g\'en\'erique $\eta$ de $\alpha^{-1}(x)$. Aussi les m\'ethodes que nous mettrons en \oe uvre viseront-elles \`a \'eliminer la ramification. Dans le cas o\`u elle est mod\'er\'ee, le lemme d'Abhyankar nous montre qu'il est possible d'y parvenir, apr\`es un nombre fini d'extensions de Kummer sur la base. Pour traiter le cas g\'en\'eral, nous utiliserons le th\'eor\`eme que d\'emontre H. Epp dans \cite{Epp}. Rappelons-en l'\'enonc\'e, sous la forme corrig\'ee qu'en proposent J. Oesterl\'e et L. Pharamond dit d'Costa (\cite{Zbarre}, appendice, th\'eor\`eme 2) :
\begin{thm}[Epp]
Soient $A$ et $A'$ deux anneaux de valuation discr\`ete, $K$ et $K'$ leur corps de fractions, $k$ et $k'$ leur corps r\'esiduel. On suppose que $A'$ domine $A$. Si la caract\'eristique $p$ de $k$ n'est pas nulle, on suppose que les \'el\'ements de $k'^{p^\infty}$, le plus grand sous-corps parfait de $k'$, sont alg\'ebriques et s\'eparables sur $k$. Il existe alors une extension alg\'ebrique $K_1$ de degr\'e fini de $K$ telle que :
\begin{enumerate}[a)]
\item la fermeture int\'egrale $A_1$ de $A$ dans $K_1$ soit un $A$-module de type fini et un anneau de valuation discr\`ete ;
\item si $K'_1$ est une extension compos\'ee de $K_1$ et $K'$, tout anneau de valuation discr\`ete $A'_1$ de corps des fractions $K'_1$ qui domine $A'$ est faiblement non ramifi\'e au-dessus de~$A_1$.
\end{enumerate}
\end{thm}
Nous aurons besoin d'utiliser le fait que la propri\'et\'e d'\^etre faiblement non ramifi\'e reste stable par certaines op\'erations. Le r\'esultat suivant se d\'eduit sans peine de la proposition $1$ de l'appendice du m\^eme article \cite{Zbarre}.  

\begin{prop}\label{fnr}
Soient $A$ et $A'$ deux anneaux de valuation discr\`ete, $K$ et $K'$ leur corps de fractions, $k$ et $k'$ leur corps r\'esiduel. On suppose que $A'$ domine $A$ et que $A'$ est faiblement non ramifi\'e au-dessus de $A$. Soit $A_1$ un anneau de valuation discr\`ete dont le corps des fractions $K_1$ est une extension alg\'ebrique de degr\'e fini de $K$ et dont le corps r\'esiduel $k_1$ est une extension s\'eparable de $k$. Alors, si $K'_1$ est une extension compos\'ee de $K_1$ et $K'$, tout anneau de valuation discr\`ete $A'_1$ de corps des fractions $K'_1$ qui domine~$A'$ est faiblement non ramifi\'e au-dessus de $A_1$.
\end{prop} 

Notons $\eta_{1},\ldots,\eta_{p}$, avec $p\in\N^*$, les points g\'en\'eriques de la fibre du morphisme $\alpha$ au-dessus de $0$. Pour chaque $i\in\cn{1}{p}$, nous allons appliquer le th\'eor\`eme de Epp aux anneaux de valuation discr\`ete $\Os_{D,0}$ et $\Os_{F,\eta_{i}}$. Les hypoth\`eses en sont v\'erifi\'ees, en vertu du r\'esultat suivant.

\begin{lem}
Supposons que la caract\'eristique du corps alg\'ebriquement clos $k$ ne soit pas nulle. Soit $\Bs$ une alg\`ebre strictement $k$-affino\"ide int\`egre. Alors le plus grand sous-corps parfait contenu dans le corps des fractions de $\Bs$ est \'egal \`a $k$.  
\end{lem}
\begin{proof}  
D'apr\`es le lemme de normalisation de Noether (\cite{BGR}, 6.1.2/2), il existe $d\in\N$ et un morphisme fini $\varphi : k\{T_1,\ldots,T_d\}\to \Bs$. La conclusion du lemme est v\'erifi\'ee pour le corps Frac($k\{T_1,\ldots,T_d\}$), car l'anneau $k\{T_1,\ldots,T_d\}$ est factoriel. D'apr\`es \cite{Epp}, \S 0.4, elle l'est encore apr\`es toute extension finie, ce qui s'applique, en particulier, \`a Frac($\Bs$).
\end{proof}

Quel que soit $i\in\cn{1}{p}$, nous obtenons ainsi une extension alg\'ebrique finie $K_{\eta_{i}}$ de Frac($k\{U\}$) v\'erifiant les conclusions du th\'eor\`eme de Epp. 

Soient $K_1$ une extension finie de $k\{U\}$ dans laquelle s'injectent tous les corps $K_{\eta_{i}}$, avec $i\in\cn{1}{p}$, et engendr\'ee par les images de ces corps. Notons $\mathscr{N}$ la fermeture int\'egrale de $k\{U\}$ dans $K_1$. Nous allons, \`a pr\'esent, consid\'erer le spectre de la fermeture int\'egrale de $\As\{T,U\}/(fT-U)$ dans un compos\'e du corps de ses fractions et de $K_1$. G\'eom\'etriquement, cela revient \`a consid\'erer le produit fibr\'e de $F$ par $\Delta=\textrm{Spec}(\mathscr{N})$ au-dessus de $D$, puis \`a le normaliser. Notons $G=\textrm{Spec}(\mathscr{G})$ le sch\'ema ainsi obtenu.

$$\xymatrix{
F \ar_\alpha[d] & F\times_D \Delta \ar_{\mu\quad  }[l] \ar_\beta[d] & G\ar^{\gamma}[ld]\ar[l]\\
D & \Delta \ar_\lambda[l]
}$$

Commen\c cons par \'enoncer quelques remarques sur les morphismes et les espaces apparaissant dans le diagramme.
\begin{enumerate}[a)] 
\item Le morphisme surjectif $\lambda : \Delta = \textrm{Spec}(\mathscr{N}) \to D$ est plat et surjectif. Puisque l'anneau $k\{U\}$ est excellent, il est \'egalement fini. En particulier, l'anneau $\mathscr{N}$ est un anneau de Dedekind et une alg\`ebre strictement $k$-affino\"ide.
\item Le morphisme $\mu : F\times_D \Delta \to F$ est, lui aussi, fini, plat et surjectif. Le morphisme $G\to F$ est donc encore fini et l'anneau $\mathscr{G}$ est une alg\`ebre strictement $k$-affino\"ide. On en d\'eduit \'egalement que toutes les composantes connexes de $G$ se surjectent sur~$F$.
\item Le morphisme $\beta : F\times_D \Delta \to \Delta$ est plat et surjectif.
\item Le morphisme $\gamma$ est surjectif. Il est \'egalement plat, puisque, quelle que soit la composante connexe $H$ de $G$, l'anneau de Dedekind $\mathscr{N}$ s'injecte dans l'anneau de $H$, qui est int\`egre et donc sans torsion.
\end{enumerate}

\'Etablissons encore deux propri\'et\'es, moins imm\'ediates : 

\begin{lem}\label{normal}
Le sch\'ema $F\times_D\Delta$ est normal hors des fibres de $\alpha\circ\mu$ au-dessus du point~$0$. En particulier, si $x$ est un point de $\Delta\setminus\lambda^{-1}(0)$, alors la fibre de $\gamma$ au-dessus de $x$ est isomorphe \`a celle de $\beta$ au-dessus de $x$ et donc \`a celle de $\alpha$ au-dessus de $\lambda(x)$.
\end{lem}
\begin{proof}
Remarquons, tout d'abord, que les fibres du morphisme $\alpha$ autres que la fibre au-dessus de $0$ sont toutes normales. En effet, pour la fibre g\'en\'erique, c'est \'evident et cela d\'ecoule de l'interpr\'etation g\'eom\'etrique des fibres du morphisme $\tau$ pour les points ferm\'es. Puisqu'en outre, $\alpha$ est plat, le morphisme $F\setminus \alpha^{-1}(0)\to D\setminus\{0\}$ est normal, au sens de \cite{EGAIV2}, 6.8.1. Puisque $\Delta\setminus\lambda^{-1}(0)$ est un sch\'ema normal, on en d\'eduit que $(F\setminus \alpha^{-1}(0))\times_{D\setminus\{0\}} (\Delta\setminus\lambda^{-1}(0))$ est encore normal, en vertu de \cite{EGAIV2}, 6.14.1. 
\end{proof}

\begin{lem}
Les fibres du morphisme $\gamma$ sont g\'eom\'etriquement r\'eduites et d\'eploy\'ees.
\end{lem} 
\begin{proof}
Les fibres du morphisme $\alpha$, \`a l'exception \'eventuelle de $\alpha^{-1}(0)$, sont g\'eom\'etriquement r\'eduites et d\'eploy\'ees. On en d\'eduit, \`a l'aide du lemme pr\'ec\'edent, que les fibres du morphisme $\gamma$ au-dessus des points de $\Delta\setminus\lambda^{-1}(0)$ le sont encore. Il nous reste \`a consid\'erer les fibres au-dessus de $\lambda^{-1}(0)$, qui est une r\'eunion finie de points ferm\'es. Puisque le corps de base $k$ est suppos\'e alg\'ebriquement clos, elles sont \'evidemment encore d\'eploy\'ees.\\ 
\indent Soient $x$ un point ferm\'e de $\Delta$ et $H$ une composante connexe de $G$. La fibre du morphisme $\gamma_H$, induit par $\gamma$ sur $H$, au-dessus de ce point est une hypersurface principale diff\'erente de $H$. Soit $\zeta$ l'un de ses points g\'en\'eriques. Il est de codimension 1 dans $H$, tout comme l'est son image $\eta$ dans $F$, par les th\'eor\`emes de Cohen et Seidenberg. Le point $\eta$ est un donc un point g\'en\'erique de la fibre de $\alpha$ au-dessus du point ferm\'e $\lambda(x)$ de $D$ et l'anneau local $\Os_{H,\zeta}$ est un anneau de valuation discr\`ete dominant $\Os_{F,\eta}$ et dont le corps des fractions co\"incide avec le corps des fonctions de $H$. La construction de $G$ et la proposition \ref{fnr} nous permettent alors d'affirmer que l'anneau de valuation discr\`ete $\Os_{G,\zeta}=\Os_{H,\zeta}$ est faiblement non ramifi\'e au-dessus de $\Os_{\Delta,x}$. Par cons\'equent, la fibre $\gamma^{-1}(x)$ est g\'en\'eriquement r\'eduite, et donc r\'eduite, puisqu'il s'agit d'une hypersurface principale d'un sch\'ema normal et noeth\'erien.
\end{proof}

Revenons, \`a pr\'esent, \`a des morphismes entre espaces analytiques. 
$$\xymatrix{
X & \Ms(\As\{T,U\}/(fT-U)) \ar[l] \ar[d]_{\tau}& \Ms(\mathscr{G})\ar[l]\ar[d]_\sigma\\
&\D & N\ar[l]_\delta
}$$
Soit $\omega$ un point rigide de $N=\Ms(\mathscr{N})$ qui s'envoie sur $0\in\D$ par le morphisme $\delta : N\to\D$. D'apr\`es le lemme \ref{normal}, il existe un voisinage affino\"ide $V$ de $\omega$ dans $N$ tel que, pour tout point rigide $v$ de $V\setminus\{\omega\}$, la fibre du morphisme $\sigma : \Ms(\mathscr{G}) \to N$ au-dessus de $v$ soit isomorphe \`a $V_\varepsilon\, \hat{\otimes}_k\, \Hs(w)$, o\`u $w =\delta(v)$, $\varepsilon = |U(w)|>0$ et $V_\varepsilon = \{x\in X\,|\, |f(x)|\ge\varepsilon\}$.\\ 
\indent Toutes les conditions sont, \`a pr\'esent, r\'eunies pour que nous puissions appliquer la th\'eor\`eme \ref{fibre} au morphisme $\sigma$ au voisinage du point $\omega$. En effet, le morphisme $\sigma$ est plat et \`a fibres g\'eom\'etriquement r\'eduites et d\'eploy\'ees, car $\gamma$ l'est. Le th\'eor\`eme nous assure l'existence d'une partie simple $P$ d'un domaine affino\"ide de $V$, de morphisme quasi-\'etales et de sections satisfaisant certaines conditions. Rappelons qu'une partie simple contenant un point rigide contient toujours un voisinage de ce point. Consid\'erons un morphisme quasi-\'etale $e : U\to V$ dont l'image contient $P$. Choisissons un point rigide $\omega'$ de $U$ qui soit un ant\'ec\'edent de $\omega$ par $e$. Puisque le corps de base $k$ est alg\' ebriquement clos, le corps r\'esiduel compl\'et\'e $\Hs(\omega)$ l'est \' egalement et le morphisme $e$ induit un isomorphisme 
$$\Hs(\omega)\xrightarrow[]{\sim} \Hs(\omega').$$
Puisque les points $\omega$ et $\omega'$ sont rigides et donc int\' erieurs, on en d\' eduit qu'il existe un isomorphisme local entre un voisinage affino\"ide de $\omega'$ et un voisinage affino\"ide de $\omega$, en vertu de \cite{bleu}, 3.4.1.\\ 
\indent En composant les diff\'erents isomorphismes r\'eciproques par les sections du th\'eor\`eme et en restreignant de fa\c con ad\'equate, nous d\'eduisons finalement l'existence d'un voisinage $V'$ de $\omega$ dans $V$ et d'une famille finie $\mathcal{T}$ de sections de $\sigma$ sur $V'$ satisfaisant la condition suivante : pour tout point $v$ de $V'$, chaque composante connexe de la fibre du morphisme $\sigma$ au-dessus de $v$ contient un et un seul \'el\'ement de la forme $t(v)$, avec $t\in\mathcal{T}$. Remarquons que nous pouvons supposer que $V'$ est connexe par arcs, puisque $N$ l'est localement.\\
\indent Puisque le morphisme $\delta : N\to \D$ est fini et plat, il est ouvert au voisinage des points rigides de $N$, donc il existe $\varepsilon>0$ tel que l'image de $V'$ par $\delta$ contienne l'ensemble $\{d\in\D\,|\, |U(d)|\le\varepsilon\}$. Soit $v$ un \'el\'ement de $V'$ tel que $|U(\delta(v))|=\varepsilon$. Soit $\varepsilon'\in\, ]0,\varepsilon]$. Il existe un chemin continu $l$ dans $V'$ joignant $v$ \`a un point $v'$ v\'erifiant les deux conditions 
$$|U(\delta(l))|\subset[\varepsilon',\varepsilon]\ \textrm{ et }\ |U(\delta(v'))|=\varepsilon'.$$
Pour chaque $t\in\mathcal{T}$, l'image du chemin $l$ par la section $t$ fournit un chemin $l_t$ dans $\Ms(\mathscr{G})$.\\
\indent Projetons, \`a pr\'esent, les chemins ainsi construits dans $X$ par le morphisme $$\Ms(\mathscr{G})\to\Ms(\As\{T,U\}/(fT-U))\to X = \Ms(\As).$$ Quel que soit $w\in l$, chaque composante connexe de $\sigma^{-1}(w)$ coupe un et un seul des chemins $l_t$, avec $t\in\mathcal{T}$. Les images de ces chemins joignent donc les composantes connexes de $V_\varepsilon$ aux composantes connexes de $V_{\varepsilon'}$. En particulier, les composantes connexes de $V_{\varepsilon}$ sont les traces de celles de $V_{\varepsilon'}$ et donc les traces de celles de $\{x\in X\, /\, f(x)\ne 0\}$. Or, d'apr\`es \cite{Bart} ou \cite{Lu}, le compl\'ementaire de l'hypersurface d\'efinie par $f$ dans $X$ est connexe, d\`es que $X$ est connexe. Par cons\'equent, le domaine affino\"ide $V_\varepsilon$ est connexe. Nous avons finalement d\'emontr\'e le r\'esultat suivant :
\begin{thm}\label{connexebis}
Soient $k$ un corps ultram\'etrique complet al\-g\'e\-bri\-que\-ment clos dont la valuation n'est pas triviale, $X$ un espace strictement $k$-affino\"ide int\`egre et $f$ une fonction analytique sur $X$ dont la norme spectrale vaut $1$. Alors le domaine affino\"ide de $X$ d\'efini par $$\{x\in X\,|\, |f(x)|\ge\varepsilon\}$$ 
est irr\'eductible, d\`es que $\varepsilon$ est assez petit. 
\end{thm}

\section{D\'emonstration des th\'eor\`emes annonc\'es}

Dans cette partie, nous expliquons comment d\'eduire les th\'eor\`emes \ref{connexe} et \ref{partition} en toute g\'en\'eralit\'e \`a partir de ceux d\'emontr\'es dans les deux paragraphes pr\'ec\'edents. Fixons $k$ un corps ultram\'etrique complet et $X$ un espace $k$-affino\"ide d'alg\`ebre $\As$.\\ 
\indent Int\'eressons-nous, tout d'abord, au th\'eor\`eme \ref{connexe}. Soient $n\in\N^*$ et $f_{1},\ldots,f_{n}$ des fonctions analytiques sur $X$. Quel que soit $\be = (\varepsilon_{1},\ldots,\varepsilon_{n}) \in (\R_{+}^*)^n$, nous noterons $V_{\be}$ le domaine analytique de $X$ d\'efini par 
$$V_{\be}= \bigcup_{1\le j\le n}\{x\in X\,|\, |f_{j}(x)|\ge\varepsilon_{j}\}.$$
Supposons que l'espace $X$ soit irr\'eductible et montrons que le domaine affino\"ide $V_{\be}$ est irr\'eductible, d\`es que $\be$ est assez petit. Comme dans le paragraphe pr\'ec\'edent, puisque le morphisme de normalisation est continu et surjectif, nous pouvons, quitte \`a remplacer $X$ par son normalis\'e, supposer que l'espace $X$ est normal.\\
\indent Soit $j\in\cn{1}{n}$. Commen\c{c}ons par nous int\'eresser aux espaces affino\"{\i}des du type 
$$V_{j,\varepsilon}=\{x\in X\,|\, |f_{j}(x)|\ge\varepsilon\},$$
avec $\varepsilon>0$. Soit $K$ un corps ultram\'etrique complet alg\'ebriquement clos contenant $k$ tel que l'espace $X\hat{\otimes}K$ soit strictement $K$-affino\"{\i}de. Le r\'esultat que nous cherchons \`a d\'emontrer est \'evident lorsque la fonction $f_{j}$ est nulle. Nous excluons dor\'enavant ce cas. D'apr\`es \cite{BGR}, 6.2.1/4 (ii), il existe alors $c\in K^*$ et $m\in\N^*$ tels que $|c\, f_{j}^m|_{\textrm{sup}}=1$. Nous pouvons donc supposer que la norme spectrale de $f_{j}$ vaut $1$, quitte \`a remplacer $f_{j}$ par $c\, f_{j}^m$, les domaines affino\"ides en jeu \'etant alors li\'es par la relation 
$$\{x\in X\hat{\otimes}_k K\,|\, |c\, f_{j}^m(x)|\ge \varepsilon\} = 
\left\{x\in X\hat{\otimes}_k K\,|\, |f_{j}(x)|\ge \left(\frac{\varepsilon}{|c|}\right)^{1/m} \right\}.$$
\indent L'espace strictement $K$-affino\"ide $X\hat{\otimes}_k K$ poss\`ede un nombre fini $Z_1,\ldots,Z_r$, avec $r\in\N$, de composantes irr\'eductibles. Sur chacune d'elles, le th\'eor\`eme \ref{connexe} est valable, d'apr\`es le th\'eor\`eme \ref{connexebis}. Par cons\'equent, il existe $\varepsilon'>0$ tel que, quel que soit $i\in\cn{1}{r}$ et quel que soit $\varepsilon\in\,]0,\varepsilon']$, l'espace 
$$\{z\in Z_i\,|\, |f_{j}(z)|\ge \varepsilon\}$$  
soit connexe.\\
\indent Pour $i\in\cn{1}{r}$, notons $Y_i$ l'image de $Z_i$ dans $X$. Quitte \`a imposer un nouvel ordre sur les indices, nous pouvons supposer qu'il existe $s\in\cn{1}{r}$ tel que $Z_1,\ldots,Z_s$ ne soient pas contenus dans 
$$V_+ = \{x\in X\,|\, f_{j}(x)\ne 0\}$$ et que $Z_{s+1},\ldots,Z_r$ le soient. Pour $i\in\cn{1}{s}$, choisissons un point $P_i$ de $Y_i$ en lequel $f_{j}$ ne s'annule pas. D'apr\`es \cite{Bart} ou \cite{Lu}, l'espace $V_+$ est connexe et m\^eme connexe par arcs, en vertu de \cite{rouge}, 3.2.1. Par cons\'equent, quel que soit $i\in\cn{2}{s}$, il existe un chemin joignant $P_1$ \`a $P_i$ sur lequel $f_{j}$ ne s'annule jamais. Ce chemin \'etant compact, la fonction $f_{j}$ y atteint son minimum $\varepsilon'_i>0$.\\
\indent Posons $\varepsilon_{j} = \min(\varepsilon',\varepsilon'_{2},\ldots,\varepsilon'_{s})$. Soit $\varepsilon\in\,]0,\varepsilon_{j}[\,$. Puisque le morphisme de changement de base $$X\hat{\otimes}_k K \to X$$ est continu et surjectif, le domaine affino\"ide $V_{j,\varepsilon}$ est connexe. D'apr\`es \cite{D}, 4.18, il est \'egalement normal, car $X$ est normal. On en d\'eduit qu'il est donc irr\'eductible.\\ 
\indent Nous pouvons, sans perte de g\'en\'eralit\'e, supposer qu'aucune des fonctions $f_{j}$, avec $j\in\cn{1}{n}$, n'est nulle. Puisque $X$ est irr\'eductible, il existe un point $x$ de $X$ en lequel aucune des fonctions $f_{j}$, avec $j\in\cn{1}{n}$, ne s'annule. Soit $\be\in \prod_{j=1}^n \,]0,\min(\varepsilon_{j},|f_{j}(x)|)[\,$. Le domaine analytique $V_{\be}$ est alors r\'eunion de parties connexes dont l'intersection contient un voisinage de $x$ dans $X$. Le lemme suivant nous montre qu'il est irr\'eductible.

\begin{lem}
Soient $V$ et $W$ deux domaines analytiques irr\'eductibles de $X$. Si l'int\'erieur de l'intersection $V\cap W$ n'est pas vide, alors la r\'eunion $V\cup W$ est irr\'eductible.
\end{lem}
\begin{proof}
Supposons que le domaine analytique $V\cup W$ soit connexe. Soient $Y$ et $Z$ deux ferm\'es de Zariski de $V\cup W$ dont la r\'eunion recouvre $V\cup W$. Supposons que $Y\ne V\cup W$. Nous avons alors $Y\cap V\ne V$ ou $Y\cap W\ne W$. Nous pouvons supposer que $Y\cap V\ne V$. Par irr\'eductibilit\'e de $V$, nous avons alors $Z\cap V=V$, autrement dit, $V\subset Z$. Par irr\'eductibilit\'e de $W$, nous devons avoir $W\subset Y$ ou $W\subset Z$. 

Supposons, par l'absurde que l'on ait $W\subset Y$. Nous avons alors $V\cap W \subset Y\cap Z$. Le domaine analytique d'int\'erieur non vide $V\cap W$ de $V$ est donc contenu dans le ferm\'e de Zariski non trivial $Y\cap Z$ du domaine analytique irr\'eductible $V$. D'apr\`es \cite{rouge}, 3.3.21, cette situation est impossible.

Finalement, nous avons $W\subset Z$ et donc $V\cup W \subset Z$. Par cons\'equent, le domaine analytique $V\cup W$ est irr\'eductible.
\end{proof}

Passons \`a la d\'emonstration du th\'eor\`eme \ref{partition}. Soit $f$ une fonction analytique sur $X$. Nous noterons $R_{X}$ le sous-$\Q$-espace vectoriel de $\R_{+}^*$ engendr\'e par les valeurs non nulles de la norme spectrale sur l'alg\`ebre $k$-affino\"{\i}de $\As$. En particulier, si $X$ est strictement $k$-affino\"{\i}de, on a $R_{X}=\sqrt{|k^*|}$. Cette d\'efinition est justifi\'ee par le lemme suivant.
\begin{lem}
Il existe un corps $L$ ultram\'etrique complet alg\'ebriquement clos et de valuation non triviale contenant $k$ tel que l'espace $X\hat{\otimes}_{k}L$ soit strictement $L$-affino\"{\i}de. Un tel corps peut \^etre choisi de fa\c{c}on \`a v\'erifier en outre 
$$|L^*|=\sqrt{|L^*|}=R_{X}.$$
\end{lem}

Dans un premier temps, nous nous int\'eresserons \`a la partie du th\'eor\`eme \ref{partition} concernant les composantes connexes. L'espace $X$ n'est plus suppos\'e irr\'eductible. Nous utiliserons la d\'efinition suivante.
\begin{defi}
Soient $Y$ un espace $k$-analytique et $g$ une fonction analytique sur $Y$. Nous dirons qu'un intervalle $I$ de $\R^+$ est \textbf{r\'egulier} pour la fonction $g$ sur l'espace $Y$ si, quels que soient $\varepsilon', \varepsilon \in I$, avec $\varepsilon'\le \varepsilon$, l'application
$$ \pi_0(\{y\in Y\,|\, |g(y)|\ge \varepsilon\}) \to \pi_0(\{y\in Y\,|\, |g(y)|\ge \varepsilon'\})$$
induite par l'inclusion est bijective.
\end{defi}
D'apr\`es \cite{D}, 5.5, le nombre de composantes connexes g\'eom\'etriques de $X$ et de ses domaines affino\"ides reste inchang\'e lorsque l'on \'etend le corps de base. Par cons\'equent, quitte \`a changer $k$ en le corps $L$ du lemme pr\'ec\'edent, nous pouvons supposer que le corps $k$ est al\-g\'e\-bri\-que\-ment clos, de valuation non triviale et que l'espace $X$ est strictement $k$-affino\"ide. Il nous faudra cependant remplacer $\sqrt{|k^*|}$ par $\sqrt{|L^*|}=|L^*|=R_{X}$.

Afin de r\'eduire encore notre probl\`eme, nous aurons besoin du lemme suivant.
\begin{lem}
Supposons que l'espace $X$ soit r\'eunion de deux ferm\'es de Zariski $Y$ et $Z$ sur lesquels il existe une partition finie de $\R^+$ en intervalles r\'eguliers pour $f$ et dont les bornes sont des \'el\'ements de $\sqrt{|k^*|}\cup\{0,+\infty\}$. Alors, la m\^eme propri\'et\'e vaut sur $X$.   
\end{lem}
\begin{proof}
Quel que soit $\varepsilon>0$, nous noterons
$$V'_\varepsilon = \{x\in Y\,|\, |f(x)|\ge \varepsilon\} \quad \textrm{ et }\quad V''_\varepsilon = \{x\in Z\,|\, |f(x)|\ge \varepsilon\}.$$
Soit $I$ un intervalle de $\R^+$ qui soit r\'egulier pour $f$ \`a la fois sur $Y$ et sur $Z$ et dont les bornes sont des \'el\'ements de $\sqrt{|k^*|}\cup\{0,+\infty\}$. Il suffit de montrer qu'un tel intervalle admet une partition finie en intervalles r\'eguliers pour $f$ sur~$X$ avec la m\^eme condition sur les bornes. Nous pouvons supposer que, quel que soit $\varepsilon\in I$, l'espace affino\"ide $V_\varepsilon$ n'est pas vide.\\
\indent Soit $\alpha\in I$. Notons $C_1,\ldots,C_r$, avec $r\in\N$, les composantes connexes de $V'_\alpha$ et $C_{r+1},\ldots,C_s$, avec $s\in\N$, celles de $V''_\alpha$. Pour $\varepsilon\in I$ et $i\in\cn{1}{r}$, nous noterons $C_{i,\varepsilon}$ l'unique composante connexe de $V'_\varepsilon$ qui v\'erifie
$$C_{i,\varepsilon} \cap V'_{\max(\varepsilon,\alpha)} = C_i \cap V'_{\max(\varepsilon,\alpha)}.$$  
Pour $\varepsilon\in I$ et $j\in\cn{r+1}{s}$, on d\'efinit de m\^eme une composante connexe $C_{j,\varepsilon}$ de $V''_\varepsilon$.\\ 
\indent Soit $\varepsilon\in I$. Remarquons que toute composante connexe $C$ de $V_\varepsilon$ s'\'ecrit de mani\`ere unique sous la forme
$$C= \bigcup_{i\in P} C_{i,\varepsilon},$$ 
o\`u $P$ d\'esigne une partie de $\cn{1}{s}$. D\'efinissons l'application 
$$c_\varepsilon : \cn{1}{s}  \to  \Pc(\cn{1}{s})$$
qui \`a un entier $i\in\cn{1}{s}$ associe l'ensemble des entiers $j\in\cn{1}{s}$ tels que $C_{j,\varepsilon}$ et $C_{i,\varepsilon}$ soient contenus dans la m\^eme composante connexe de $V_\varepsilon$.\\
\indent L'application $c : \varepsilon \mapsto c_\varepsilon$ ne peut prendre qu'un nombre fini de valeurs et est d\'ecroissante, au sens o\`u, pour $\varepsilon'\ge \varepsilon$, on a
$$\forall i\in\cn{1}{s},\, c_{\varepsilon'}(i)\subset c_\varepsilon(i).$$
Par cons\'equent, il existe une partition finie $\Pc$ de $I$ en intervalles sur lesquels l'application $c$ est constante. Chacun de ces intervalles est r\'egulier pour $f$ sur $X$.\\
\indent Soit $\beta \in\sqrt{|k^*|}$ tel que l'application $c$ soit constante sur l'intervalle $I\cap\mathopen[0,\beta\mathclose]$. Les parties $C_{1,\beta},\ldots,C_{s,\beta}$ sont alors des domaines strictement affino\"ides de $V_\beta$ et donc de $X$. Nous pouvons choisir les intervalles de la partition $\Pc$ de fa\c con que leurs bornes diff\'erentes de celles de l'intervalle $I$ soient contenues dans l'ensemble $E$ des \'el\'ements $\varepsilon$ de $I$ pour lesquels il existe des indices $i,j\in\cn{1}{s}$ tels que
$$\forall \varepsilon'\in I\cap [0,\varepsilon\mathclose[,\, C_{i,\beta}\cap C_{j,\beta}\cap V_{\varepsilon'}\ne \emptyset$$
et $$\forall \varepsilon'\in I\cap \mathopen]\varepsilon,+\infty\mathclose[,\, C_{i,\beta}\cap C_{j,\beta}\cap V_{\varepsilon'} = \emptyset.$$
En d'autres termes, chaque \'el\'ement de $E$ peut \^etre obtenu comme la valeur maximale de la valeur absolue de la fonction $f$ sur un certain domaine strictement affino\"ide de $X$. D'apr\`es \cite{BGR}, 6.2.1/4, on en d\'eduit que $E\subset \sqrt{|k^*|}$. 
\end{proof}

T\^achons, tout d'abord, de d\'emontrer qu'il existe une partition finie de $\R^+$ en intervalles r\'eguliers pour $f$ sur $X$ et dont les bornes sont des \'el\'ements de $\sqrt{|k^*|}\cup\{0,+\infty\}$. Puisque l'espace $k$-affino\"ide $X$ poss\`ede un nombre fini de composantes irr\'eductibles, le lemme pr\'ec\'edent nous montre qu'il suffit de le prouver pour chacune d'elles. Nous pouvons donc supposer que l'espace $X$ est irr\'eductible et m\^eme int\`egre. D'apr\`es \cite{BGR}, 6.2.1/4 (ii), si la fonction $f$ n'est pas nulle, nous pouvons supposer que sa norme spectrale vaut~$1$. Dans ce cas, nous savons, d'apr\`es le th\'eor\`eme \ref{connexe}, qu'il existe $\varepsilon'\in\,]0,1]\cap \sqrt{|k^*|}$ tel que, quel que soit $\varepsilon\in\,]0,\varepsilon']$, l'espace 
$$V_\varepsilon = \{x\in X\,|\, |f(x)|\ge \varepsilon\}$$  
soit connexe. Par hypoth\`ese, l'espace $V_0=X$ est connexe, donc l'intervalle $[0,\varepsilon'[$ est r\'egulier pour $f$ sur $X$. Quel que soit $\varepsilon \in\, ]1,+\infty[$, l'espace $V_\varepsilon$ est vide et l'intervalle $]1,+\infty[$ est donc \'egalement r\'egulier pour $f$ sur $X$. Le th\'eor\`eme \ref{variation} nous assure encore qu'il est possible de d\'ecouper l'intervalle $[\varepsilon',1]$ en un nombre fini d'intervalles r\'eguliers pour $f$ sur $X$ et dont les bornes sont des \'el\'ements de $\sqrt{|k^*|}\cup\{0,+\infty\}$. Par cons\'equent, il existe une partition finie de $\R^+$ en intervalles r\'eguliers pour $f$ sur $X$ dont les bornes jouissent de la m\^eme propri\'et\'e.\\
\indent Int\'eressons-nous \`a pr\'esent \`a la forme des intervalles de la partition pr\'ec\'edente. Les lemmes qui suivent nous permettront de l'obtenir, concluant ainsi la d\'emonstration du th\'eor\`eme \ref{partition}.
\begin{lem} \label{injection}
Quel que soit $\beta>0$, il existe $\alpha\in[0,\beta[$ tel que, quel que soit $\varepsilon \in[\alpha,\beta]$, l'application naturelle 
$$\pi_0(V_\beta) \to \pi_0(V_\varepsilon)$$
soit injective.
\end{lem}
\begin{proof}
Soit $\beta>0$. Puisque $V_\beta$ ne poss\`ede qu'un nombre fini de composantes connexes, il suffit de montrer que deux d'entre elles, $C_0$ et $C_1$, distinctes, sont contenues dans deux composantes connexes distinctes de $V_\varepsilon$, avec $\varepsilon\le\beta$, d\`es que $\varepsilon$ est assez proche de $\beta$.\\
\indent Remarquons qu'il existe une fonction, $g$ analytique sur $V_\beta$, v\'erifiant 
$$g_{|C_0}\equiv 0,\ g_{|C_1}\equiv 1 \textrm{ et } g^2-g=0.$$
Puisque $V_\beta$ est un domaine rationnel de $X$, nous pouvons approcher la fonction $g$, uniform\'ement sur $V_\beta$, par une suite de quotients d'\'el\'ements de $\As$ sans p\^oles sur $V_\beta$. Par cons\'equent, il existe $p,q\in\As$ tels que $q$ ne s'annule pas sur $V_\beta$ et 
$$\forall x\in V_\beta,\, |g(x)-h(x)|\le \frac{1}{3} \textrm{ et } |h^2(x)-h(x)|\le \frac{1}{5},$$
o\`u $h=p/q$.
Le lieu d'annulation de $q$ est une partie compacte de $X$, disjointe de $V_\beta$, sur laquelle la fonction continue $f$ atteint son maximum $M<\beta$. Soit $M'\in\mathopen]M,\beta]$. La fonction m\'eromorphe $h$ est analytique sur $V_{M'}$. D\'efinissons un compact $K$ de $V_{M'}$ par 
$$K = \left\{x\in V_{M'}\,|\, |h^2(x)-h(x)|\ge \frac{1}{4}\right\}.$$
La fonction continue $f$ y atteint son maximum $M_1$. Puisque $K$ et $V_\beta$ sont disjoints, on a n\'ecessairement $M_1<\beta$. Fixons $\alpha\in\mathopen]M_1,\beta\mathclose[$.\\
\indent Soit $\varepsilon \in[\alpha,\beta]$. Le compact $K$ est disjoint de $V_\varepsilon$ donc, quel que soit $x\in V_\varepsilon$, on a $|h(h-1)(x)|< 1/4$. Posons 
$$D_0=\left\{x\in V_\varepsilon\,|\, |h(x)|<\frac{1}{2}\right\} \textrm{ et } D_1=\left\{x\in V_\varepsilon\,|\, |h(x)-1|<\frac{1}{2}\right\}.$$
Ces deux ouverts sont disjoints et recouvrent $V_\varepsilon$. Par cons\'equent, ils sont r\'eunions de composantes connexes. En outre, $C_0\subset D_0$ et $C_1 \subset D_1$, donc les parties $C_0$ et $C_1$ sont contenues dans deux composantes connexes distinctes de $V_\varepsilon$.
\end{proof}
\begin{lem}
Si l'intervalle $\mathopen[\alpha,\beta\mathclose[$ est r\'egulier pour $f$ sur $X$, alors l'intervalle $\mathopen[\alpha,\beta\mathclose]$ l'est encore.
\end{lem}
\begin{proof}
Il nous suffit de montrer que l'application naturelle
$$\iota : \pi_0(V_\beta) \to \pi_0(V_\alpha)$$
est bijective. D'apr\`es le lemme \ref{injection}, elle est injective. Montrons qu'elle est \'egalement surjective.\\
\indent Soit $C$ une composante connexe de $V_\alpha$. C'est une partie compacte sur laquelle la fonction continue $f$ atteint son maximum $M$. Puisque l'intervalle $\mathopen[\alpha,\beta\mathclose[$ est r\'egulier pour $f$ sur $X$, on a $C\cap V_\varepsilon \ne\emptyset$ et donc $M\ge\varepsilon$, quel que soit $\varepsilon \in \mathopen[\alpha,\beta\mathclose[$. On en d\'eduit que $M\ge\beta$, autrement dit que $C\cap V_\beta \ne\emptyset$. Choisissons une composante connexe de $V_\beta$ coupant $C$. Elle s'envoie sur $C$ par $\iota$. L'application $\iota$ est donc surjective.
\end{proof}
\indent Il nous reste \`a d\'emontrer la partie du th\'eor\`eme \ref{partition} qui concerne les composantes irr\'eductibles. Elle se d\'eduit de celle qui concerne les composantes connexes lorsqu'on l'applique au normalis\'e de $X$.\\ 

\'Enon\c cons, \`a pr\'esent, un corollaire du th\'eor\`eme \ref{connexe}. Il figure d\'ej\`a dans \cite{rouge}, \S 2.3, sans d\'emonstration.

\begin{cor}
Un point d'un bon espace $k$-analytique en lequel l'anneau local est int\`egre poss\`ede une base de voisinages affino\"ides irr\'eductibles. 
\end{cor}
\begin{proof}
L'espace \'etant bon, il suffit de d\'emontrer que tout point d'un espace $k$-affino\"ide en lequel l'anneau local est int\`egre poss\`ede un voisinage affino\"ide irr\'eductible. Soient $Y$ un espace $k$-affino\"ide et $y$ un point de $Y$ en lequel l'anneau local $\Os_{Y,y}$ est int\`egre. Alors le point $y$ ne peut \^etre situ\'e que sur une seule des composantes irr\'eductibles de $Y$. Notons $F$ cette composante et $G$ la r\'eunion des autres. Soit $g$ une fonction analytique sur $Y$ nulle sur le ferm\'e de Zariski $G$ et ne s'annulant pas en $y$. D'apr\`es le th\'eor\`eme \ref{connexe}, il existe $\varepsilon\in\, ]0,|f(x)|[$ tel que le domaine affino\"ide
$$\{y'\in Y\,|\, |f(y')|\ge \varepsilon\} = \{y'\in F\,|\, |f(y')|\ge\varepsilon\}$$
soit irr\'eductible.  
\end{proof}

\section{Privil\`ege}

Dans cette partie, nous \'enon\c cons et d\'emontrons un r\'esultat de privil\`ege pour les vari\'et\'es analytiques $p$-adiques. La septi\`eme partie de l'article \cite{Douady} d'A. Douady est consacr\'ee \`a cette notion, dans le cadre analytique complexe. Rappelons-en quelques d\'efinitions et notations.\\
\indent Si $K$ est une partie compacte de $\C^n$, avec $n\in\N$, on note $\Os(K)$ l'espace vectoriel des germes de fonctions analytiques au voisinage de $K$ et $B(K)$ son adh\'erence dans l'espace de Banach des fonctions continues sur $K$. Si $\Fs$ est un faisceau analytique coh\'erent d\'efini au voisinage de $K$, on note $\Fs(K)$ la limite inductive des modules des sections de $\Fs$ sur les voisinages ouverts de $K$ et
$$B(K,\Fs) = B(K)\otimes_{\Os(K)} \Fs(K).$$
\indent Revenons, \`a pr\'esent, au cadre des espaces analytiques d\'efinis sur un corps ultram\'etrique complet. Soient $Y$ un espace $k$-analytique normal et s\'epar\'e et $V=\Ms(\As)$ un domaine affino\"ide de $Y$ contenu dans l'int\'erieur de $Y$. En d\'efinissant $B(V)$ de la m\^eme fa\c con que pr\'ec\'edemment, on obtient un isomorphisme
$$B(V) \xrightarrow{\sim} \As$$ 
et l'on retrouve les sections du faisceau structural pour la G-topologie. De m\^eme, si $\Fs$ d\'esigne un faisceau coh\'erent pour la G-topologie de $Y$, la formule d\'efinissant $B(V,\Fs)$ redonne exactement le $\As$-module de type fini $\Fs(V)$ des sections globales de $\Fs$ sur $V$.\\

\begin{defi}
Soient $k$ un corps ultram\'etrique complet, $Y$ un bon espace $k$-analytique et $\Fs$ un faisceau coh\'erent d\'efini sur $Y$. Nous dirons qu'un voisinage affino\"ide $V$ d'un point $y$ de l'espace $k$-analytique $Y$ est \textbf{privil\'egi\'e} pour le faisceau $\Fs$ s'il v\'erifie 
$$\Fs(V) \hookrightarrow \Fs_y,$$
o\`u $\Fs(V)$ doit \^etre pris au sens de la G-topologie et $\Fs_y$ au sens de la topologie sur l'espace topologique sous-jacent $|Y|$.  
\end{defi}

Pour qu'un voisinage compact $K$ d'un point soit privil\'egi\'e, A. Douady impose non seulement la condition qui figure dans la d\'efinition, mais encore une autre qui porte sur des propri\'et\'es d'exactitude du foncteur $B(K,.)$ (\emph{cf.} \cite{Douady}, \S 7, d\'efinition 2). Dans notre cadre, elles seront toujours v\'erifi\'ees pour les domaines affino\"ides.\\

Signalons que l'on peut penser \`a un voisinage privil\'egi\'e pour un faisceau coh\'erent comme un voisinage sur lequel vaut une g\'en\'eralisation de l'unicit\'e du prolongement analytique. En effet, si $Y$ d\'esigne un espace $k$-analytique irr\'eductible et r\'eduit, nous savons, d'apr\`es \cite{rouge}, 3.3.21, qu'une fonction nulle sur un ouvert non vide de $Y$ est identiquement nulle. Ce r\'esultat se traduit par le fait que $Y$ est un voisinage privil\'egi\'e de tous ses points pour le faisceau structural. Nous en d\'eduisons aussit\^ot le lemme suivant.

\begin{lem} \label{irr}
Dans un espace $k$-analytique r\'eduit, un voisinage affino\"ide d'un point est privil\'egi\'e pour le faisceau structural d\`es que toutes les composantes irr\'eductibles du voisinage passent par ce point.
\flushright{$\square$}
\end{lem}

\indent Nous souhaitons montrer ici que tout point d'un espace $k$-analytique poss\`ede un syst\`eme fondamental de voisinages affino\"ides privil\'egi\'es pour un faisceau coh\'erent fix\'e, du moins lorsque l'espace est bon. Ce r\'esultat est analogue \`a celui d\'emontr\'e par A. Douady dans \cite{Douady} (\S 6, th\'eor\`eme 1). Notre d\'emonstration reprend des id\'ees qui figurent dans l'article \cite{Frisch} de J. Frisch.\\

\begin{lem}\label{se}
Soient $Y$ un espace $k$-analytique, $y$ un point de $Y$ et $$0\to \Fs'\to\Fs\to\Fs''$$ une suite exacte de faisceaux coh\'erents sur $Y$. Alors, un voisinage affino\"ide de $y$ privil\'egi\'e pour les faisceaux $\Fs'$ et $\Fs''$ l'est encore pour le faisceau $\Fs$.
\end{lem}
\begin{proof}
Ce r\'esultat provient directement de l'exactitude \`a gauche du foncteur des sections globales.
\end{proof}

\begin{lem}\label{filtration}
Soient $Y$ un espace $k$-affino\"ide et $\Fs$ un faisceau coh\'erent sur $Y$. Alors il existe un entier $r\in\N$, une filtration
$$0=\Fs_0\subset \Fs_1\subset \ldots \subset \Fs_r=\Fs$$
de $\Fs$ par des sous-faisceaux coh\'erents et $r$ ferm\'es de Zariski de $Y$ int\`egres $Z_0,\ldots,Z_{r-1}$ v\'erifiant la condition suivante : quel que soit $i\in\cn{0}{r-1}$, on dispose d'un isomorphisme de faisceaux
$$\Fs_{i+1}/\Fs_i \simeq \Os_{Z_i}.$$
\end{lem}
\begin{proof}
Le module $\Fs(Y)$ des sections du faisceau coh\'erent $\Fs$ sur $Y$ est un module de type fini sur l'alg\`ebre $\Bs$ de $Y$. Par cons\'equent, il existe un entier $r\in\N$, une filtration
$$0=M_0\subset M_1\subset \ldots \subset M_r=M$$
de $M$ par des sous-$\Bs$-modules de type fini v\'erifiant la condition suivante : quel que soit $i\in\cn{0}{r-1}$, il existe un id\'eal premier $\p_i$ de $\Bs$ et un isomorphisme 
$$M_{i+1}/M_i \simeq \Bs/\p_i.$$
Pour $i\in\cn{0}{r}$, notons $\Fs_i$ le faisceau coh\'erent associ\'e \`a $M_i$ sur $Y$ et, pour $i\in\cn{0}{r-1}$, notons $Z_i$ le ferm\'e de Zariski int\`egre de $Y$ d'alg\`ebre $\Bs/\p_i$. Ils satisfont la conclusion du lemme.
\end{proof}

\begin{thm}
Soient $k$ un corps ultram\'etrique complet, $Y$ un bon espace $k$-analytique et $\Fc$ une famille finie de faisceaux coh\'erents sur $Y$. Tout point de $Y$ poss\`ede un syst\`eme fondamental de voisinages affino\"ides privil\'egi\'es pour chacun des faisceaux de $\Fc$. 
\end{thm}
\begin{proof}
Soit $y\in Y$. Par d\'efinition d'un bon espace, le point $y$ poss\`ede un syst\`eme fondamental $\Vs$ de voisinages affino\"ides dans $Y$. Il suffit de montrer que, quel que soit $V\in\Vs$, le point $y$ poss\`ede un voisinage affino\"ide dans $V$ qui soit privil\'egi\'e pour chacun des faisceaux de $\Fc$. Soit $V\in\Vs$. Notons $\Bs$ son alg\`ebre.\\
\indent Soit $\Fs$ un \'el\'ement de $\Fc$. D'apr\`es le lemme \ref{filtration}, il existe entier $r(\Fs)\in\N$, une filtration
$$0=\Fs_0\subset \Fs_1\subset \ldots \subset \Fs_{r(\Fs)}=\Fs$$
de $\Fs$ par des sous-faisceaux coh\'erents et $r(\Fs)$ ferm\'es de Zariski de $V$ int\`egres $$Z_{\Fs,0},\ldots,Z_{\Fs,r(\Fs)-1}$$ v\'erifiant la condition suivante : quel que soit $i\in\cn{0}{r(\Fs)-1}$, on dispose d'un isomorphisme de faisceaux
$$\Fs_{i+1}/\Fs_i \simeq \Os_{Z_{\Fs,i}}.$$
\indent D\'efinissons l'ensemble
$$\Pc = \{Z_{\Fs,i},\, \Fs\in\Fc, 0\le i\le r(\Fs)-1\}.$$
D'apr\`es le lemme \ref{se}, il nous suffit pour conclure de montrer que le point $y$ poss\`ede un voisinage affino\"ide privil\'egi\'e pour chacun des faisceaux $\Os_Z$, avec $Z\in \Pc$. Notons $\Qc$ l'ensemble des \'el\'ements de $\Pc$ \'evitant le point $y$. Leur r\'eunion $R$ d\'efinit un ferm\'e de Zariski de $V$ ne contenant pas $y$. Par cons\'equent, il existe une fonction $g\in \Bs$ qui soit nulle sur $R$, mais pas en $y$. D'apr\`es le th\'eor\`eme \ref{connexe}, il existe $\varepsilon \in\, ]0,|g(y)|[$ tel que l'espace affino\"ide $$\{z\in Z\,|\, |g(z)|\ge \varepsilon \}$$ soit irr\'eductible, quel que soit $Z$ dans $\Pc\setminus \Qc$. Notons $W$ le voisinage affino\"ide de $y$ dans $V$ d\'efini par
$$W=\{z\in V\,|\, |g(z)|\ge \varepsilon\}.$$
Soit $Z\in\Pc$. Si $Z\in \Qc$, le faisceau $\Os_Z$ restreint \`a $W$ est nul et le voisinage $W$ de $y$ est donc privil\'egi\'e pour $\Os_Z$. Si $Z\notin \Qc$, le ferm\'e de Zariski $Z\cap W$ est irr\'eductible et on dispose donc d'un morphisme injectif
$$\Os_Z(W) \simeq \Os_{Z\cap W}(Z\cap W) \hookrightarrow \Os_{Z,y}.$$
Autrement dit, le voisinage $W$ de $y$ est, dans ce cas encore, privil\'egi\'e pour $\Os_Z$.
\end{proof}

\appendix

\section{Un analogue $p$-adique du th\'eor\`eme de J. Frisch}

Nous proposons ici un analogue, dans le cadre des espaces analytiques d\'efinis sur un corps ultram\'etrique complet, du th\'eor\`eme I,9 qui figure dans l'article \cite{Frisch} de J. Frisch. Nous suivrons, ici, la d\'emonstration de C. B\u anic\u a et O. St\u an\u a\c sil\u a (\emph{cf.} \cite{BS2}, 5, fin du \S 3). Nous obtiendrons une version un peu plus g\'en\'erale du th\'eor\`eme, proche de celle que propose Y.-T. Siu dans \cite{Siu}. Commen\c cons par une d\'efinition et un lemme.

\begin{defi}
Soient $k$ un corps ultram\'etrique complet et $Y$ un espace $k$-analytique. Une partie $A$ de $Y$ est dite \textbf{morcelable} si, pour tout ferm\'e de Zariski $Z$ d\'efini au voisinage de $A$, l'image r\'eciproque de $A \cap Z$ dans le normalis\'e de $Z$ poss\`ede un nombre fini de composantes connexes.
\end{defi}

\begin{lem}\label{support}
Soient $k$ un corps ultram\'etrique complet et $Y$ un espace $k$-analytique. Si l'espace $Y$ est normal, alors le support de tout faisceau d'id\'eaux coh\'erent sur $Y$ est ouvert et ferm\'e.
\end{lem}
\begin{proof}
Nous pouvons supposer que $Y$ est un espace $k$-affino\"ide. Soit $\Js$ un faisceau d'id\'eaux coh\'erents sur $Y$. Comme tout faisceau coh\'erent, son support est ferm\'e dans $Y$. Pour conclure, il nous suffit de montrer qu'il est \'egalement ouvert. Soit $y$ un point de $Y$ en lequel la fibre de $\Js$ n'est pas nulle. Alors il existe un voisinage ouvert connexe $V$ de $y$ et une fonction analytique $g\in \Js(V)$ qui ne soit pas identiquement nulle sur $V$. Soit $z\in V$. D'apr\`es \cite{rouge}, 3.3.21, le ferm\'e de Zariski d\'efini par $g$ est d'int\'erieur vide dans l'espace normal et connexe $V$. En particulier, la fonction $g$ n'est nulle sur aucun voisinage de $z$ dans $V$. On en d\'eduit que la fibre $\Js_z$ n'est pas nulle et donc que le support de $\Js$ est ouvert.    
\end{proof}

Le r\'esultat de finitude sur lequel nous nous appuierons concerne les familles croissantes de faisceaux coh\'erents.

\begin{lem}
Soient $k$ un corps ultram\'etrique complet et $Y$ un espace $k$-analytique. Notons $\pi : \tilde{Y}\to Y$ le morphisme de normalisation. Soit $A$ une partie de l'espace topologique sous-jacent \`a $Y$ telle que $\pi^{-1}(A)$ poss\`ede un nombre fini de composantes connexes. Soit $(\Is_n)_{n\in\N}$ une suite croissante de faisceaux d'id\'eaux coh\'erents de $\Os_Y$ d\'efinis chacun sur un voisinage de $A$. Soit $a$ un point de $A$ en lequel la fibre $(\Is_n)_a$ est nulle, quel que soit $n\in\N$. Alors il existe un voisinage de $a$ dans $A$ sur lequel toutes les fibres du faisceau $\Is_n$ sont nulles, quel que soit $n\in\N$. 
\end{lem}
\begin{proof}
Quel que soit $b\in \pi^{-1}(a)$, il existe un voisinage $V_b$ de $b$ dans $\tilde{Y}$ tel que $\pi^{-1}(A)\cap V_b$ soit connexe. La partie 
$V=\bigcup_{b\in\pi^{-1}(A)} V_b$
d\'efinit un voisinage de la fibre $\pi^{-1}(a)$ dans $\tilde{Y}$, donc il existe un voisinage $U$ de $a$ dans $Y$ tel que $\pi^{-1}(U)\subset V$. Nous allons montrer que, quel que soit $n\in\N$ et quel que soit $a'\in U\cap A$, on a 
$$(\Is_n)_{a'}=0.$$
Soit $n\in\N$. Notons $\Js_n$ le faisceau d\'efini par 
$$\Js_n = \pi^{-1}(\Is_n)\, \Os_{\tilde{Y}}.$$
C'est un faisceau d'id\'eaux coh\'erent d\'efini sur un voisinage de $\pi^{-1}(A)$ dans $\tilde{Y}$. Quel que soit $b\in \pi^{-1}(a)$, il existe un voisinage $V_{b,n}$ de $\pi^{-1}(A)\cap V_b$ sur lequel $\Js_n$ est d\'efini. Puisque 
$$(\Js_n)_b = (\Is_n)_a\, \Os_{\tilde{Y},b}=0,$$
la fibre de $\Js_n$ est nulle en tout point de $V_{b,n}$, d'apr\`es le lemme \ref{support}. La partie 
$V_n = \bigcup_{b\in\pi^{-1}(A)} V_{b,n}$
est un voisinage de $\pi^{-1}(A\cap U)$ dans $\tilde{Y}$. Par cons\'equent, la fibre du faisceau $\pi_*\Js_n$ est nulle en tout point de $A\cap U$. Or le diagramme commutatif
$$\xymatrix{
\Js_n \ar[d]\ar[r]& \pi_*\Js_n \ar[d]\\
\Os_Y \ar[r] & \pi_*\Os_{\tilde{Y}}
}$$
montre que le faisceau $\Js_n$ s'injecte dans le faisceau $\pi_*\Js_n$. Le r\'esultat annonc\'e s'en d\'eduit.
\end{proof}

\begin{prop}
Soient $k$ un corps ultram\'etrique complet, $Y$ un espace $k$-affino\"ide et $A$ une partie de l'espace topologique sous-jacent \`a $Y$. Soient $\Fs$ un faisceau coh\'erent d\'efini sur $Y$ et $(\Fs_n)_{n\in\N}$ une suite croissante de sous-faisceaux coh\'erents de $\Fs$ d\'efinis chacun sur un voisinage affino\"ide de $A$. Si la partie $A$ est morcelable, alors la suite $(\Fs_n)_{n\in\N}$ est localement stationnaire dans $A$ au sens o\`u, quel que soit $a\in A$, il existe un entier $n_0\in\N$ et un voisinage $U$ de $a$ dans $A$ tels que 
$$\forall n\ge n_0,\, \forall z\in A,\, (\Fs_{n_0})_z\xrightarrow{\sim} (\Fs_n)_z.$$
\end{prop}
\begin{proof}
Supposons que la partie $A$ soit morcelable. Soit $a\in A$. Il existe $n_0\in \N$ tel que, quel que soit $n\ge n_0$, on ait 
$$(\Fs_{n_0})_a\xrightarrow{\sim} (\Fs_n)_a.$$
Quitte \`a restreindre $Y$, \`a remplacer $\Fs$ par $\Fs/\Fs_{n_0}$ et $\Fs_n$ par $\Fs_n/\Fs_{n_0}$, pour $n\ge n_0$, puis \`a d\'ecaler les indices, nous pouvons supposer que $$(\Fs_n)_a = 0,$$ quel que soit $n\in\N$. D'apr\`es le lemme \ref{filtration}, il existe un entier $r\in\N$, une filtration
$$0=\Fs^{(0)}\subset \Fs^{(1)}\subset \ldots \subset \Fs^{(r)}=\Fs$$
de $\Fs$ par des sous-faisceaux coh\'erents et $r$ ferm\'es de Zariski de $Y$ int\`egres $Z_0,\ldots,Z_{r-1}$ v\'erifiant la condition suivante : quel que soit $i\in\cn{0}{r-1}$, on dispose d'un isomorphisme de faisceaux
$$\Fs^{(i+1)}/\Fs^{(i)} \simeq \Os_{Z_i}.$$
Il nous suffit, \`a pr\'esent, de montrer que, pour chaque $i\in\cn{0}{r-1}$, la sous-suite $(\Gs_{i,n})_{n\in\N}$ de $\Fs^{(i)}/\Fs^{(i+1)} \simeq \Os_{Z_i}$ induite par $(\Fs_n)_{n\in\N}$ stationne au voisinage de $a$ dans $A$ et m\^eme au voisinage de $a$ dans $A \cap Z_i$. Le lemme pr\'ec\'edent nous permet de conclure.
\end{proof}

Il ne nous reste plus qu'\`a rendre global le r\'esultat pr\'ec\'edent pour obtenir le th\'eor\`eme recherch\'e. 

\begin{thm}
Soient $k$ un corps ultram\'etrique complet, $Y$ un bon espace $k$-analytique et $K$ une partie compacte de l'espace topologique sous-jacent \`a $Y$. Si $K$ est morcelable et poss\`ede un syst\`eme fondamental de voisinages affino\"ides, alors l'anneau $\Os(Y,K)$ des germes de fonctions analytiques au voisinage de $K$ est noeth\'erien.
\end{thm}
\begin{proof}
Soit $(I_n)_{n\in\N}$ une suite croissante d'id\'eaux de type fini de $\Os(Y,K)$. Pour $n\in\N$, notons $\Is_n$ le faisceau d'id\'eaux coh\'erents de $\Os_Y$ engendr\'e par $I_n$. D'apr\`es la proposition pr\'ec\'edente, la suite $(\Is_n)_{n\in\N}$ stationne sur $K$, au sens o\`u il existe $n_0\in\N$ tel que, quel que soit $n\ge n_0$ et quel que soit $y\in K$, on dispose d'un isomorphisme 
$$(\Is_{n_0})_y \xrightarrow{\sim} (\Is_n)_y.$$
\indent Puisque l'id\'eal $I_{n_0}$ est fini, il poss\`ede un syst\`eme g\'en\'erateur fini $(f_1,\ldots,f_p)$, avec $p\in\N$ et $f_i\in \Os(Y,K)$, quel que soit $i\in\cn{1}{p}$. Le morphisme de faisceaux 
$$\varphi : \begin{array}{ccc}
\Os_Y^p & \to & \Is_{n_0}\\
(a_1,\ldots,a_p) & \mapsto & a_1 f_1 + \ldots + a_p f_p
\end{array}$$
est surjectif.\\
\indent Soit $n\ge n_0$. Si un morphisme entre deux faisceaux coh\'erents induit un isomorphisme entre les fibres en un point, alors il induit un isomorphisme au voisinage de ce point. Par cons\'equent, les faisceaux $\Is_{n_0}$ et $\Is_n$ co\"incident sur un voisinage $U_n$ de $K$. Soient $g\in I_n$ et $V_n$ un voisinage affino\"ide de $K$ dans $U_n$ sur lequel les fonctions $f_1,\ldots, f_p,g$ soient d\'efinies. Notons $\Gs$ le noyau du morphisme de faisceaux $\varphi$. C'est encore un faisceau coh\'erent sur $V_n$. De la suite exacte $0\to \Gs \to \Os^p \to \Is_n \to 0$, on d\'eduit une surjection 
$$\begin{array}{ccc}
\Os(V_n)^p & \to & \Is_n(V_n)\\
(a_1,\ldots,a_p) & \mapsto & a_1 f_1 + \ldots + a_p f_p
\end{array},$$ 
car $H^1(V_n,\Gs)=0$. Par cons\'equent, 
$$g\in (f_1,\ldots,f_p)\, \Os(V_n) \subset  (f_1,\ldots,f_p)\, \Os(Y,K) = I_{n_0}.$$
On en d\'eduit que $I_n=I_{n_0}$.
\end{proof}

Signalons que le th\'eor\`eme que d\'emontre J. Frisch concerne des compacts poss\'edant un syst\`eme fondamental de voisinages compos\'e d'espaces de Stein. Dans le cadre des espaces analytiques d\'efinis sur un corps ultram\'etrique complet $k$, il existe \'egalement une notion d'espace de Stein. Un espace $k$-analytique $Y$ est dit de Stein s'il existe une suite croissante $(Y_n)_{n\in\N}$ de domaines affino\"ides de $Y$ v\'erifiant les conditions suivantes :
\begin{enumerate}[a)]
\item quel que soit $n\in\N$, $Y_n$ est un domaine de {Weierstra\ss} de $Y_{n+1}$ ;
\item la famille $\{Y_n,\, n\in\N\}$ d\'efinit un G-recouvrement de $Y$.
\end{enumerate}
La seconde condition, pr\'esente dans \cite{AB}, 2.3, fait d\'efaut dans \cite{rouge}, p. 96.\\
\indent Il d\'ecoule de la d\'efinition qu'une partie compacte d'un espace analytique poss\'edant un syst\`eme fondamental de voisinages constitu\'e d'espaces de Stein poss\`ede encore un syst\`eme fondamental de voisinages constitu\'e d'affino\"ides.\\

\indent Mentionnons, pour conclure, deux exemples de parties compactes morcelables : 
\begin{enumerate}[a)]
\item Nous dirons qu'une partie $A$ d'un espace analytique $Y$ est semi-analytique si tout point de $A$ poss\`ede un voisinage affino\"ide dans lequel la partie $A$ est semi-alg\'ebrique, c'est-\`a-dire d\'ecrite par un nombre fini d'in\'egalit\'es entre fonctions. Une telle partie est morcelable, lorsqu'elle est compacte, d'apr\`es \cite{salg}, 3.2. Dans ce cas, nous retrouvons exactement l'\'enonc\'e original de J. Frisch.
\item Si $Y$ est un espace $k$-analytique et $K$ une extension de $k$, nous noterons $Y_K$ l'espace $K$-analytique obtenu par extension du corps de base. Nous dirons qu'un morphisme $\varphi$ entre espaces $k$-analytiques est une immersion s'il se d\'ecompose sous la forme
$$\varphi : Z \hookrightarrow Y_K \to Y,$$
o\`u $Z$ d\'esigne un ferm\'e de Zariski d'un domaine analytique de $Y_K$ et s'il induit un hom\'eomorphisme de $Z$ sur son image $Z'$ et des isomorphismes entre les corps r\'esiduels compl\'et\'es en tous les points de $Z'$.
L'image de toute immersion d\'efinit une partie morcelable. Les fibres des morphismes entre espaces $k$-analytiques rentrent, par exemple, dans ce cadre. Remarquons qu'elles ne sont pas semi-analytiques, en g\'en\'eral.
\end{enumerate}

\nocite{*}
\bibliographystyle{plain}
\bibliography{biblio}


\end{document}